\documentclass[preprint,a4paper]{elsarticle}

\usepackage[lmargin=27mm, tmargin=36mm, body={156mm, 225mm}]{geometry}

\usepackage{amsmath,amsthm,amssymb,mathrsfs}

\newtheorem{theorem}{Theorem}[section]
\newtheorem{corollary}[theorem]{Corollary}
\newtheorem{remark}[theorem]{Remark}
\newtheorem{lemma}[theorem]{Lemma}

\newtheorem{example}[theorem]{Example}

\numberwithin{equation}{section}

\journal{}

\allowdisplaybreaks

\begin{document}

\title{Sphere theorems for submanifolds in  K\"ahler Manifold\tnoteref{SS}}

\author[whu]{Jun Sun}
\ead{sunjun@whu.edu.cn}

\author[whu]{Linlin Sun\corref{sll}}
\ead{sunll@whu.edu.cn}

 \tnotetext[SS]{The first author was supported by the National Natural Science Foundation of China (Grant No. 11401440). Part of the work was finished when the first author was a visiting scholar at MIT supported by China Scholarship Council (CSC) and Wuhan University. The author would like to express his gratitude to Professor Tobias Colding for his invitation, to MIT for their hospitality, and to CSC and Wuhan University for their support. The second author was supported by the National Natural Science Foundation of China (Grant No. 11801420)  and Fundamental Research Funds for the Central Universities (Grant No. 2042018kf0044). }
 
 \address[whu]{School of Mathematics and Statistics \& Computational Science Hubei Key Laboratory, Wuhan University, 430072 Wuhan, China}
 
 \cortext[sll]{Corresponding author.}
 
\begin{abstract}

In this paper, we prove some differentiable sphere theorems and topological sphere theorems for submanifolds in K\"ahler manifold, especially in complex space forms.

\end{abstract}


\maketitle
\section{Introduction}

\noindent The study of the relation between curvature and topology is a fundamental problem in differential geometry. Sphere theorems play an important role in such a study.  There are two types of differentiable sphere theorems: one is for the Riemannian manifold itself (i.e., intrinsic version), the other is for submanifolds in a Riemannian manifold (i.e., extrinsic version). The typical example of the former one is the classical $1/4$-pinched differentiable sphere theorem, which states that a compact Riemannian manifold $M$ of dimension $n\geq 4$ with pointwise $1/4$-pinched sectional curvature is diffeomorphic to a spherical space form. This theorem was finally proved by Brendle-Schoen (\cite{BS}, \cite{BS2}).

We are mainly interested in the latter one, i.e., extrinsic version. There are many interesting sphere theorems for smooth submanifold $M$ immersed into a Riemannian manifold $N^{n+p}$. For example, Lawson-Simons (\cite{LS}) considered the vanishing theorem of integral current in an $n$-dimensional submanifold in unit sphere (the case of submanifold in Euclidean space was considered by Xin (\cite{Xin}) and showed that an $n$-dimensional submanifold in unit sphere with $|{\bf B}|^2<\min\{n-1,2\sqrt{n-1}\}$ is a homotopy sphere. Shiohama and Xu (\cite{SXu}) improved Lawson-Simons' result to complete submanifold in space forms with nonnegative sectional curvature. Cui-Sun (\cite{CS}) and Gu-Xu (\cite{GX}, \cite{XG1}, \cite{XG2}, etc.) also proved some topological and differentiable sphere theorems for submanifolds in general Riemannian manifold. Furthermore, Li-Wang (\cite{LW}) proved some differentiable sphere theorems for Lagrangian submanifolds in complex space form. In general, the conditions of sphere theorems are expressed in terms of the scalar curvature, Ricci curvature  or the sectional curvature and the mean curvature of the submanifold and the sectional curvature of the ambient manifold.

\vspace{.1in}

In this paper, we will consider  sphere theorems for submanifolds in K\"ahler manifold. Contrary to the above mentioned sphere theorems, we will express the condition in terms of the holomorphic sectional curvature of the ambient manifold instead of its sectional curvature.

\vspace{.1in}

Let $M$ be a smooth $n$-dimensional submanifold of a  K\"ahler manifold $N^{2m}$. We will denote the curvature tensors on $M$ and $N$ by $R$ and $K$, respectively. Recall that the sectional curvature is given by
\begin{equation*}
K(X,Y):=K(X,Y,X,Y)
\end{equation*}
and the holomorphic sectional curvature is given by
\begin{equation*}
K(X):=K(X,JX):=K(X,JX,X,JX),
\end{equation*}
where $X$ and $Y$ are tangent vector fields on $M$. Denote the minimal and maximal holomorphic sectional curvatures by
\begin{equation}\label{e1.4}
\tilde{K}_{\min}:=\min_{|X|=1}K(X),  \quad \tilde{K}_{\max}:=\max_{|X|=1}K(X).
\end{equation}
Our first theorem is as follows:

\vspace{.1in}

\noindent \textbf{Theorem A:} {\it Let $M$ be a smooth $n(\geq 2)$-dimensional closed simply connected submanifold of a  K\"ahler manifold $N^{2m}$. If the scalar curvature of $M$ satisfies the following condition:
\begin{equation}\label{e-A}
R_M\geq
      \begin{cases}
            \ \ \ \ \frac{3n^2+8}{4}\tilde{K}_{\max}-\frac{n^2-n+4}{2}\tilde{K}_{\min}+\frac{n-2}{n-1}|{\bf H}|^2, \ \ if \ \tilde{K}_{\min}\geq 0;\\
            \ \ \ \ \frac{3n^2+8}{4}\tilde{K}_{\max}-\frac{n^2-n+8}{2}\tilde{K}_{\min}+\frac{n-2}{n-1}|{\bf H}|^2, \ \ if \ \tilde{K}_{\min}\leq 0\leq \tilde{K}_{\max};\\
            \frac{3(n^2-n+2)}{4}\tilde{K}_{\max}-\frac{n^2-n+8}{2}\tilde{K}_{\min}+\frac{n-2}{n-1}|{\bf H}|^2, \ if \ \tilde{K}_{\max}\leq 0,
      \end{cases}
\end{equation}
and we further assume that the strict inequality holds for some point $x_0\in M$ if $\tilde{K}_{\max}=\tilde{K}_{\min}$. Then $M$ is diffeomorphic to ${\mathbb S}^n$.
}

\vspace{.1in}

Recall that a submanifold $M$ in a K\"ahler manifold $N$ is said to be {\it totally real} in $N$ if $JT_x(M)\subset N_x(M)$ for each $x\in M$, where $J$ is the complex structure on $N$ and $N_x(M)$ is the normal space of $M$ in $N$ at $x$. When the submanifold is totally real, Theorem A can be improved to be the following:

\vspace{.1in}

\begin{corollary}\label{cor-totallyreal}
Let $M$ be a smooth $n(\geq 2)$-dimensional closed simply connected totally real submanifold of a K\"ahler manifold $N^{2m}$. If $M$ satisfies the following condition:
\begin{equation*}
R_M\geq\frac{3(n^2-n+2)}{4}\tilde{K}_{\max}-\frac{n^2-n+4}{2}\tilde{K}_{\min}+\frac{n-2}{n-1}|{\bf H}|^2,
\end{equation*}
and we further assume that the strict inequality holds for some point $x_0\in M$ if $\tilde{K}_{\max}=\tilde{K}_{\min}=0$. Then $M$ is diffeomorphic to ${\mathbb S}^n$.
\end{corollary}

\vspace{.1in}

In particular, when $N$ is a complex space form with constant holomorphic sectional curvature $c$, we have:

\vspace{.1in}

\begin{corollary}\label{cor1.2}
Let $M$ be a smooth $n(\geq 2)$-dimensional closed simply connected totally real submanifold of complex space form $N^{2m}$ with holomorphic sectional curvature $c$. If $M$ satisfies  the following condition:
\begin{equation*}
R_M\geq\frac{(n-2)(n+1)}{4}c+\frac{n-2}{n-1}|{\bf H}|^2,
\end{equation*}
and we further assume that the strict inequality holds for some point $x_0\in M$ if $c=0$. Then $M$ is diffeomorphic to ${\mathbb S}^n$.
\end{corollary}

\vspace{.1in}

Next, we plan to examine differentiable sphere theorems under Ricci curvature pinching condition. 

\vspace{.1in}

\noindent \textbf{Theorem B:} {\it For fixed $0<\varepsilon\leq1$, set $\delta(\varepsilon,n)=\frac{\left((n-4)\varepsilon+2\right)^2}{4\left(2+\left(n^2-4n+2\right)\varepsilon\right)}$. Let $M$ be a smooth $n(\geq 4)$-dimensional closed simply connected submanifold of  a  K\"ahler manifold $N^{2m}$. If $M$ satisfies the following condition:
\begin{equation}\label{e-B}
Ric^{[2]}_{\min}\geq
      \begin{cases}
          \ \ \  \frac{3n+4\varepsilon}{2}\tilde{K}_{\max}-(n-1+2\varepsilon)\tilde{K}_{\min}+\delta(\varepsilon,n)|{\bf H}|^2, \ \ if \ \tilde{K}_{\min}\geq 0;\\
         \ \  \ \frac{3n+4\varepsilon}{2}\tilde{K}_{\max}-(n-1+4\varepsilon)\tilde{K}_{\min}+\delta(\varepsilon,n)|{\bf H}|^2, \ \ if \ \tilde{K}_{\min}\leq 0\leq \tilde{K}_{\max};\\
           \frac{3(n-1+\varepsilon)}{2}\tilde{K}_{\max}-(n-1+4\varepsilon)\tilde{K}_{\min}+\delta(\varepsilon,n)|{\bf H}|^2, \ \ if \ \tilde{K}_{\max}\leq 0,
      \end{cases}
\end{equation}
and the strict inequality holds for some point $x_0\in M$. Then $M$ is diffeomorphic to ${\mathbb S}^n$.
}

\vspace{.1in}

\begin{corollary}\label{cor1.7}
For fixed $0<\varepsilon\leq1$, set $\delta(\varepsilon,n)=\frac{\left((n-4)\varepsilon+2\right)^2}{4\left(2+\left(n^2-4n+2\right)\varepsilon\right)}$. Let $M$ be a smooth $n(\geq 4)$-dimensional closed simply connected totally real submanifold of a K\"ahler manifold $N^{2m}$. If $M$ satisfies  the following condition:
\begin{equation*}
Ric^{[2]}_{\min}\geq\frac{3(n-1+\varepsilon)}{2}\tilde{K}_{\max}-(n-1+2\varepsilon)\tilde{K}_{\min}+\delta(\varepsilon,n)|{\bf H}|^2,
\end{equation*}
and the strict inequality holds for some point $x_0\in M$. Then $M$ is diffeomorphic to ${\mathbb S}^n$.
\end{corollary}

\begin{corollary}\label{cor1.8}
For fixed $0<\varepsilon\leq1$, set $\delta(\varepsilon,n)=\frac{\left((n-4)\varepsilon+2\right)^2}{4\left(2+\left(n^2-4n+2\right)\varepsilon\right)}$. Let $M$ be a smooth $n(\geq 4)$-dimensional closed simply connected totally real submanifold of  complex space form $N^{2m}$ with holomorphic sectional curvature $c$. If $M$ satisfies  the following condition:
\begin{equation*}
Ric^{[2]}_{\min}\geq \frac{n-1-\varepsilon}{2}c+\delta(\varepsilon,n)|{\bf H}|^2,
\end{equation*}
and the strict inequality holds for some point $x_0\in M$. Then $M$ is diffeomorphic to ${\mathbb S}^n$.
\end{corollary}

\begin{remark}\label{remark1.5}
If $\varepsilon=1$, then $\delta(\varepsilon, n)=\frac{1}{4}$.
\end{remark}

\vspace{.1in}

For a submanifold in a K\"ahler manifold, we also have the following topological sphere theorem:

\vspace{.1in}

\noindent \textbf{Theorem C:} {\it Let $M$ be a smooth $n(\geq 4)$-dimensional closed simply connected submanifold of a  K\"ahler manifold $N^{2m}$. If the scalar curvature of $M$ satisfies the following condition:
\begin{equation}\label{e-C}
R_M\geq
      \begin{cases}
         \ \  \ \    \frac{3n^2+16}{4}\tilde{K}_{\max}-\frac{n^2-n+8}{2}\tilde{K}_{\min}+\frac{n-3}{n-2}|{\bf H}|^2, \ \ if \ \tilde{K}_{\min}\geq 0;\\
         \ \ \ \    \frac{3n^2+16}{4}\tilde{K}_{\max}-\frac{n^2-n+16}{2}\tilde{K}_{\min}+\frac{n-3}{n-2}|{\bf H}|^2, \ \ if \ \tilde{K}_{\min}\leq 0\leq \tilde{K}_{\max};\\
            \frac{3(n^2-n+4)}{4}\tilde{K}_{\max}-\frac{n^2-n+16}{2}\tilde{K}_{\min}+\frac{n-3}{n-2}|{\bf H}|^2, \ \ if \ \tilde{K}_{\max}\leq 0,
      \end{cases}
\end{equation}
and the strict inequality holds for some point $x_0\in M$. Then $M$ is homeomorphic to ${\mathbb S}^n$.
}

\vspace{.1in}

\begin{corollary}\label{Cor1.6}
Let $M$ be a smooth $n(\geq 4)$-dimensional closed simply connected totally real submanifold of a K\"ahler manifold $N^{2m}$. If $M$ satisfies the following condition:
\begin{equation*}
R_M\geq\frac{3(n^2-n+4)}{4}\tilde{K}_{\max}-\frac{n^2-n+8}{2}\tilde{K}_{\min}+\frac{n-3}{n-2}|{\bf H}|^2,
\end{equation*}
and the strict inequality holds for some point $x_0\in M$. Then $M$ is homeomorphic to ${\mathbb S}^n$.
\end{corollary}

\vspace{.1in}

In particular, when $N$ is a complex space form with constant holomorphic sectional curvature $c$, we have (comparing with Corollary \ref{cor1.2}):

\begin{corollary}\label{Cor1.7}
Let $M$ be a smooth $n(\geq 4)$-dimensional closed simply connected totally real submanifold of complex space form $N^{2m}$ with holomorphic sectional curvature $c$. If $M$ satisfies  the following condition:
\begin{equation*}
R_M\geq\frac{n^2-n-4}{4}c+\frac{n-3}{n-2}|{\bf H}|^2,
\end{equation*}
and the strict inequality holds for some point $x_0\in M$. Then $M$ is homeomorphic to ${\mathbb S}^n$.
\end{corollary}

\vspace{.1in}

\noindent \textbf{Theorem D:} {\it Let $M$ be a smooth $n(\geq 4)$-dimensional closed simply connected submanifold of a  K\"ahler manifold $N^{2m}$. If $M$ satisfies the following condition:
\begin{equation}\label{e-D}
Ric^{[4]}_{\min}\geq
      \begin{cases}
            (3n+4)\tilde{K}_{\max}-2(n+1)\tilde{K}_{\min}+\frac{1}{2}|{\bf H}|^2,  \ \ if \ \tilde{K}_{\min}\geq 0;\\
            (3n+4)\tilde{K}_{\max}-2(n+3)\tilde{K}_{\min}+\frac{1}{2}|{\bf H}|^2, \ \ if \ \tilde{K}_{\min}\leq 0\leq \tilde{K}_{\max};\\
            \ \ \ \ \ \ \ 3n\tilde{K}_{\max}-2(n+3)\tilde{K}_{\min}+\frac{1}{2}|{\bf H}|^2, \ \ if \ \tilde{K}_{\max}\leq 0,
      \end{cases}
\end{equation}
and the strict inequality holds for some point $x_0\in M$. Then $M$ is homeomorphic to ${\mathbb S}^n$.
}

\vspace{.1in}

\begin{corollary}\label{cor1.12}
Let $M$ be a smooth $n(\geq 4)$-dimensional simply connected compact totally real submanifold of a K\"ahler manifold $N^{2m}$. If $M$ satisfies the following condition:
\begin{equation*}
Ric^{[4]}_{\min}\geq 3n\tilde{K}_{\max}-2(n+1)\tilde{K}_{\min}+\frac{1}{2}|{\bf H}|^2,
\end{equation*}
and the strict inequality holds for some point $x_0\in M$. Then $M$ is homeomorphic to ${\mathbb S}^n$.
\end{corollary}

In particular, if $N$ is a complex space form, then we have the following topological sphere theorem for totally real submanifold (comparing with Remark \ref{remark1.5}):

\begin{corollary}\label{Cor1.9}
Let $M$ be a smooth $n(\geq 4)$-dimensional simply connected compact totally real submanifold of a complex space form $N^{2m}$ with holomorphic sectional curvature $c$. If $M$ satisfies the following condition:
\begin{equation*}
Ric^{[4]}_{\min}\geq (n-2)c+\frac{1}{2}|{\bf H}|^2,
\end{equation*}
and the strict inequality holds for some point $x_0\in M$. Then $M$ is homeomorphic to ${\mathbb S}^n$.
\end{corollary}

\begin{remark}
All results mentioned above are sharp. 
\begin{itemize}
\item Consider the totally embedding $\mathbb{C}P^{n/2}(4)\subset\mathbb{C}P^m(4)$ where $n$ is an even number. Then $Ric=(n+2)g$ and $R_M=n(n+2)$. Thus Theorem A, Theorem B, Theorem C and Theorem D are sharp.
\item Consider $M_{p,\mu}:=\mathrm{S}^{n-p}\left(\frac{\mu}{\sqrt{1+\mu^2}}\right)\times\mathrm{S}^{p}\left(\frac{1}{\sqrt{1+\mu^2}}\right)\left(\subset\mathrm{S}^{n+1}(1)\right)\subset\mathbb{C}P^{n+1}(4)$ where $0<\mu<1$, then $M_{p,\mu}$ is a totally real submanifold of $\mathbb{C}P^{n+1}(4)$. Moreover,
\begin{align*}
R_{M_{1,\mu}}-\dfrac{n-2}{n-1}|{\bf H}|^2-(n-2)(n+1)=-\dfrac{n-2}{n-1}\mu^2\to0, \quad\text{as}\ \mu\to0,\\
R_{M_{2,\mu}}-\dfrac{n-3}{n-2}|{\bf H}|^2-(n^2-n-4)=-\dfrac{2(n-4)}{n-2}\mu^2\to0, \quad\text{as}\ \mu\to0.
\end{align*}
Therefore, Corollary \ref{cor-totallyreal}, Corollary \ref{cor1.2}, Corollary \ref{Cor1.6} and Corollary \ref{Cor1.7} are optimal. 
\item For $\varepsilon=1$, Corollary \ref{cor1.7}, Corollary \ref{cor1.8} are optimal for $n=4$. Corollary \ref{cor1.12} and Corollary \ref{Cor1.9} are optimal for $n=4$. We refer the reader to \cite{XG2}.

\end{itemize}
\end{remark}

\vspace{.1in}

In another paper, we will consider differentiable sphere theorems and topological sphere theorems for Lagrangian submanifods in K\"ahler manifold (\cite{SS1}).
Similar argument can also prove some sphere theorems for submanifolds in Sasaki space forms.

\vspace{.2in}

\section{Preliinaries}

\vspace{.1in}

In this section, we will provide some basic materials about K\"ahler manifold that will be used in the proof of the main theorems. First recall the following expression of the sectional curvature and curvature tensor in terms of holomorphic sectional curvature:

\begin{lemma}[cf. \cite{Kar}]
 Let $N$ be a Riemannian manifold and $X$, $Y$, $Z$, $W$ be vector fields on $N$. Then we have
\begin{eqnarray}\label{E2.1}
24K(X,Y,Z,W)
&=& \ \ K(X+Z,Y+W)+K(X-Z,Y-W)\nonumber\\
& & +K(X+W,Y-Z)+K(X-W,Y+Z)-K(X+Z,Y-W) \nonumber\\
& & -K(X-Z,Y+W)-K(X+W,Y+Z)-K(X-W,Y-Z).
\end{eqnarray}
\end{lemma}

\begin{lemma}[cf. \cite{YK}]
Let $N$ be a K\"ahler manifold and $X$, $Y$ be vector fields on $N$. Then we have
\begin{equation}\label{E2.2}
32K(X,Y)=3K(X+JY)+3K(X-JY)-K(X+Y)-K(X-Y)-4K(X)-4K(Y).
\end{equation}
\end{lemma}

Putting (\ref{E2.2}) into (\ref{E2.1}), we get that

\begin{corollary}
Let $N$ be a K\"ahler manifold and $X$, $Y$, $Z$, $W$ be vector fields on $N$. Then we have
\begin{eqnarray}\label{E2.3}
256K(X,Y,Z,W)
&=&\ \  K(X+Z+JY+JW)+K(X+Z-JY-JW)\nonumber\\
& & -K(X+Z+JY-JW)-K(X+Z-JY+JW) \nonumber\\
& & +K(X-Z+JY-JW)+K(X-Z-JY+JW)\nonumber\\
& & -K(X-Z+JY+JW)-K(X-Z-JY-JW) \nonumber\\
& & +K(X+W+JY-JZ)+K(X+W-JY+JZ)\nonumber\\
& & -K(X+W+JY+JZ)-K(X+W-JY-JZ) \nonumber\\
& & +K(X-W+JY+JZ)+K(X-W-JY-JZ)\nonumber\\
& & -K(X-W+JY-JZ)-K(X-W-JY+JZ).
\end{eqnarray}
\end{corollary}

\vspace{.1in}

Let $M^n$ be an $n$-dimensional submanifold in Riemannian manifold $N^{d}$. Choose local orthonormal frame $\{e_1,\cdots,e_{d}\}$ on $N$ so that $\{e_1,\cdots,e_n\}$ are tangent to $M$ and $\{e_{n+1},\cdots,e_{d}\}$ are normal to  $M$. Denote $R$ and $K$ the Levi-Civita connections on $M$ and $N$, respectively, and $h^{\alpha}_{ij}=\langle {\bf B}(e_i,e_j),e_{\alpha}\rangle$ the component of the second fundamental form of $M$ in $N$. The mean curvature vector is given by ${\bf H}=\sum_{\alpha=n+1}^dH^{\alpha}e_{\alpha}$, where $H^{\alpha}=\sum_{i=1}^nh^{\alpha}_{ii}$.
Then the Gauss equation can be written as
\begin{equation}\label{e-gauss}
R_{ijkl}=K_{ijkl}+\sum_{\alpha=n+1}^d(h^{\alpha}_{ik}h^{\alpha}_{jl}-h^{\alpha}_{il}h^{\alpha}_{jk}).
\end{equation}
In particular, the Ricci curvature and the scalar curvature satisfies
\begin{equation*}
Ric(e_i)=R_{ii}=\sum_{j=1}^nK_{ijij}+\sum_{\alpha=n+1}^d\sum_{j=1}^n[h^{\alpha}_{ii}h^{\alpha}_{jj}-(h^{\alpha}_{ij})^2],
\end{equation*}
\begin{equation}\label{e-scalar}
R_M=\sum_{i,j=1}^nK_{ijij}+|{\bf H}|^2-|{\bf B}|^2.
\end{equation}

\vspace{.1in}

Fix $p\in M$, $X,Y\in T_pM$ and an orthonormal basis $\{e_1,\cdots,e_n\}$ of $T_pM$, the following notations will be used in this paper:
\begin{equation*}
Ric(X,Y)=\sum_{i=1}^nR(X,e_i,Y,e_i), \ \ Ric_{jj}=Ric(e_j,e_j),
\end{equation*}
\begin{equation*}
[e_{i_1},\cdots,e_{i_k}]=span\{e_{i_1},\cdots,e_{i_k}\}, \ \ \ \forall 1\leq i_1<i_2 < \cdots <i_k\leq n,
\end{equation*}
\begin{equation*}
Ric^{[k]}[e_{i_1},\cdots,e_{i_k}]=\sum_{j=1}^kRic_{i_ji_j}, \ \ Ric^{[k]}_{\min}(p)=\min_{[e_{i_1},\cdots,e_{i_k}]\subset T_pM}Ric^{[k]}[e_{i_1},\cdots,e_{i_k}],
\end{equation*}
where $Ric^{[k]}[e_{i_1},\cdots,e_{i_k}]$ is called the {\it $k$-th weak Ricci curvature} of $[e_{i_1},\cdots,e_{i_k}]$, which was first introduced by Gu-Xu in \cite{GX}.

\vspace{.1in}

At the end of this section, we will state some lemmas which will be crucial in the proof of our main theorems. The first result is due to Aubin:

\begin{lemma}[\cite{Au}]\label{lemma-Aubin}Let $M$ be a  compact n-dimensional Riemannian manifold. If $M$ has nonnegative Ricci curvature everywhere and has positive Ricci curvature at some point, then $M$ admits a metric with positive Ricci curvature everywhere.
\end{lemma}

A Riemannian manifold $M$ is said to have {\it nonnegative (positive, respectively) isotropic curvature}, if
\begin{equation*}
R_{1313}+R_{1414}+R_{2323}+R_{2424}-2R_{1234}\geq 0 (>0, respectively)
\end{equation*}
for all orthonormal four-frames $\{e_1,e_2.e_3.e_4\}$. This conception was introduced by Micallef-Moore and they proved the following topological sphere theorem:

\begin{lemma}[\cite{MM}]\label{lemma-MM} Let $M$ be a  compact simply connected $n(\geq 4)$-dimensional Riemannian manifold which has positive isotropic curvature, then $M$ is homeomorphic to a sphere.
\end{lemma}

In addition, Micallef-Wang proved the following topological result for manifold with positive isotropic curvature:

\begin{lemma}[\cite{MW}]\label{lemma-MW} Let $M$ be a closed even-dimensional Riemannian manifold which has positive isotropic curvature, then $b_2(M)=0$.
\end{lemma}

Furthermore, Seshadri proved the following result for manifold with nonnegative isotropic curvature:

\begin{lemma}[\cite{Se}]\label{lemma-Se}Let $M$ be a compact $n$-dimensional Riemannian manifold. If $M$ has nonnegative isotropic curvature everywhere and has positive isotropic curvature at some point, then $M$ admits a metric with positive isotropic curvature.
\end{lemma}

The 1/4-differentiable sphere theorem was finally proved by Brendle-Schoen (\cite{BS}, \cite{BS2}) using the Ricci flow method. They proved that:

\begin{theorem}[\cite{BS}]\label{thmBS}
Let $(M,g_0)$ be a compact, locally irreducible Riemannian manifold of dimension $n(\geq 4)$ with curvature tensor $R$. Assume that $M\times{\mathbb R}^2$ has nonnegative isotropic curvature, i.e.,
\begin{equation}\label{e-BS}
R_{1313}+\lambda^2R_{1414}+\mu^2R_{2323}+\lambda^2\mu^2R_{2424}-2\lambda\mu R_{1234}\geq0
\end{equation}
for all orthonormal four-frames $\{e_1,e_2,e_3,e_4\}$ and all $\lambda, \mu\in[-1,1]$. Then one of the following statements holds:

\ \ (i) $M$ is diffeomorphic to a spherical space form;

\ (ii) $n=2m$ and the universal covering of $M$ is a K\"ahler manifold biholomorphic to ${\mathbb C}P^m$;

(iii) The universal covering of $M$ is isometric to a compact symmetric space.
\end{theorem}

\vspace{.2in}

\section{Some algebraic estimates }

\vspace{.1in}

In this section, we will prove some algebraic estimates that are used in the proof of the main theorems.

In this section, we always assume $n\geq4$. We say that $R$ is an algebraic curvature on $\mathbb{R}^n$ if $R$ is a fourth tensor such that for every $x,y,z,w\in\mathbb{R}^n$,
\begin{align*}
\begin{cases}
R(x,y,z,w)=-R(y,x,z,w)=-R(x,y,w,z)=R(z,w,x,y),\\
R(x,y,z,w)+R(y,z,x,w)+R(z,x,y,w)=0.
\end{cases}
\end{align*}

Let $\left\{e_i\right\}_{i=1}^n$ be an orthonormal frame of $\mathbb{R}^n$.
\begin{example}
 If $B=(h_{ij}^{\alpha}):\mathbb{R}^n\times\mathbb{R}^n\longrightarrow\mathbb{R}^p$ is a bilinear operator, we obtain an algebraic curvature tensor $\tilde R$ defined by:
\begin{align*}
\tilde R_{ijkl}:=\sum_{\alpha=1}^ph_{ik}^{\alpha}h_{jl}^{\alpha}-\sum_{\alpha=1}^ph_{il}^{\alpha}h_{jk}^{\alpha},\quad\forall 1\leq i, j, k, l\leq n.
\end{align*}
\end{example}

\begin{lemma}\label{lem31}
Let $R$ be an algebraic curvature tensor $R$. Suppose there is a constant $c$ such that for every orthonormal four-frames $\{e_1,e_2,e_3,e_4\}$,
\begin{align*}
R_{1212}+R_{1234}\geq c,
\end{align*}
then for every $\lambda, \mu\in[-1,1]$ and  every orthonormal four-frames $\{e_1,e_2,e_3,e_4\}$
\begin{align*}
R_{1313}+\lambda^2R_{1414}+\mu^2R_{2323}+\lambda^2\mu^2R_{2424}-2\lambda\mu R_{1234}\geq\left(1+\lambda^2\right)\left(1+\mu^2\right)c.
\end{align*}
\end{lemma}

\noindent\textbf{Proof:}
The assumption implies that for every orthonormal four-frames $\{e_1,e_2,e_3,e_4\}$,
\begin{equation*}
R_{1212}-|R_{1234}| \geq c.
\end{equation*}
The Bianchi identity yields that
\begin{align*}
R_{1234}=&R_{1324}+R_{1432}.
\end{align*}
Therefore,
\begin{align*}
&R_{1313}+\lambda^2R_{1414}+\mu^2R_{2323}+\lambda^2\mu^2R_{2424}-2\lambda\mu R_{1234}\\
=&R_{1313}+\lambda^2R_{1414}+\mu^2R_{2323}+\lambda^2\mu^2R_{2424}-2\lambda\mu(R_{1324}+R_{1432})\\
\geq&R_{1313}+\lambda^2R_{1414}+\mu^2R_{2323}+\lambda^2\mu^2R_{2424}-\left(1+\lambda^2\mu^2\right)|R_{1324}|-\left(\lambda^2+\mu^2\right)|R_{1432}|\\
=&\left(R_{1313}-|R_{1324}|\right)+\lambda^2\left(R_{1414}-|R_{1432}|\right)+\mu^2\left(R_{2323}-|R_{2314}|\right)+\lambda^2\mu^2\left(R_{2424}-|R_{2413}|\right)\\
\geq&\left(1+\lambda^2+\mu^2+\lambda^2\mu^2\right)c\\
=&\left(1+\lambda^2\right)\left(1+\mu^2\right)c.
\end{align*}

\hfill Q.E.D.

\begin{lemma}Let $R$ be an algebraic curvature tensor $R$. Suppose there is a constant $c$ such that for every orthonormal four-frames $\{e_1,e_2,e_3,e_4\}$,
\begin{align*}
R_{1313}+R_{2323}+R_{1234}\geq c,
\end{align*}
then for every $\lambda\in[-1,1]$ and  every orthonormal four-frames $\{e_1,e_2,e_3,e_4\}$
\begin{align*}
R_{1313}+\lambda^2R_{1414}+R_{2323}+\lambda^2R_{2424}-2\lambda R_{1234}\geq \left(1+\lambda^2\right)c.
\end{align*}
\end{lemma}

\noindent\textbf{Proof:} A straightforward verification.
\hfill Q.E.D.

\begin{lemma}\label{lem:a-2}Let $B=(h_{ij}^{\alpha}):\mathbb{R}^n\times\mathbb{R}^n\longrightarrow\mathbb{R}^p$ is a bilinear operator. Define $H^{\alpha}:=\sum_{i=1}^nh^{\alpha}_{ii}$ and
\begin{align*}
\tilde R_{ijkl}:=\sum_{\alpha=1}^ph_{ik}^{\alpha}h_{jl}^{\alpha}-\sum_{\alpha=1}^ph_{il}^{\alpha}h_{jk}^{\alpha},\quad\forall 1\leq i, j, k, l\leq n.
\end{align*}
Then  for all orthonormal four-frames $\{e_1,e_2,e_3,e_4\}$, we have
\begin{equation}\label{eA1}
\tilde R_{1212}+\tilde R_{1234}\geq \dfrac12\left[\frac{\sum_{\alpha=1}^p\left(H^{\alpha}\right)^2}{n-1}-\sum_{i,j=1}^n\sum_{\alpha=1}^p\left(h_{ij}^{\alpha}\right)^2\right],
\end{equation}
with equality holds if and only if $h^{\alpha}_{ii}=h^{\alpha}_{11}+h^{\alpha}_{22}$ for all $i\neq 1,2$ and $h^{\alpha}_{ij}=0$ for all distinct $i,j$ with $\{i,j\}\neq \{1,2\}$.
We also have
\begin{equation}\label{eA2}
\sum_{i=1}^2\sum_{j=3}^4\tilde R_{ijij}-2\tilde R_{1234}\geq \frac{\sum_{\alpha=1}^p\left(H^{\alpha}\right)^2}{n-2}-\sum_{i,j=1}^n\sum_{\alpha=1}^p\left(h_{ij}^{\alpha}\right)^2.
\end{equation}
\end{lemma}

\noindent\textbf{Proof:} For the proof of this Lemma, we refer the reader to Gu-Xu's paper \cite{GX}. We only need to notice that (\ref{eA1}) follows from the inequality
 \begin{equation}\label{e-diff}
 2h^{\alpha}_{mm}h^{\alpha}_{ll}\geq \sum_{i\neq j}(h^{\alpha}_{ij})^2+\frac{(H^{\alpha})^2}{n-1}-\sum_{i,j=1}^n(h^{\alpha}_{ij})^2,
\end{equation}
for all distinct $m,l$, and the equality holds if and only if 
 \begin{equation}\label{e-diff2}
h^{\alpha}_{ii}=h^{\alpha}_{mm}+h^{\alpha}_{ll}, \ \ for \ all \ i\neq m,l.
\end{equation}
Furthermore, (\ref{eA2}) follows from the inequality
 \begin{equation}\label{e-homm}
 2h^{\alpha}_{pp}h^{\alpha}_{qq}+2h^{\alpha}_{mm}h^{\alpha}_{ll}\geq \sum_{i\neq j}(h^{\alpha}_{ij})^2+\frac{(H^{\alpha})^2}{n-2}-\sum_{i,j=1}^n(h^{\alpha}_{ij})^2,
\end{equation}
for all distinct $p,q,m,l$, and the equality holds if and only if 
 \begin{equation}\label{e-homm2}
h^{\alpha}_{ii}=h^{\alpha}_{pp}+h^{\alpha}_{qq}=h^{\alpha}_{mm}+h^{\alpha}_{ll}, \ \ for \ all \ i\neq p,q,m,l.
\end{equation}
\hfill Q.E.D.

\begin{lemma}\label{lem:a-3}Let $B$ and $\tilde R$ be as in Lemma \ref{lem:a-2}. Assume 
\begin{align*}
\sum_{k=1}^n\tilde R_{ikik}+\sum_{k=1}^n\tilde R_{jkjk}\geq 2D,\quad\forall 1\leq i<j\leq n,
\end{align*}
then for every $0<\varepsilon\leq 1$ and all orthonormal four-frames $\{e_1,e_2,e_3,e_4\}$
\begin{align*}
\tilde R_{1212}+\tilde R_{1234}\geq\frac{1}{\varepsilon}\left[D-\frac{\left((n-4)\varepsilon+2\right)^2}{8\left(2+\left(n^2-4n+2\right)\varepsilon\right)}\sum_{\alpha=1}\left(H^{\alpha}\right)^2\right].
\end{align*}
\end{lemma}

\noindent\textbf{Proof:} The proof can be found in \cite{CS}. For reader's convenience, we give another but direct proof.
Set
\begin{align*}
h^{\alpha}_{ij}:=\mathring{h}^{\alpha}_{ij}+\frac{1}{n}H^{\alpha}\delta_{ij},\quad T^{\alpha}:=\dfrac{1}{n}H^{\alpha}.
\end{align*}
One can check that
\begin{align*}
\sum_{i,j=1}^n\left(h_{ij}^{\alpha}\right)^2=\sum_{i,j=1}^n\left(\mathring{h}_{ij}^{\alpha}\right)^2+n\left(T^{\alpha}\right)^2,\quad\forall 1\leq\alpha\leq p.
\end{align*}
Denoted by $\tilde R_{ii}:=\sum_{j=1}^n\tilde R_{ijij}$, we get
\begin{align*}
\tilde R_{ii}=&H^{\alpha}h^{\alpha}_{ii}-\sum_{j=1}^n\sum_{\alpha=1}^ph^{\alpha}_{ij}h^{\alpha}_{ij}\\
=&(n-1)\sum_{\alpha=1}^p\left(T^{\alpha}\right)^2+\sum_{\alpha=1}^p\left[(n-2)T^{\alpha}\mathring{h}^{\alpha}_{ii}-\sum_{j=1}^n\left(\mathring{h}^{\alpha}_{ij}\right)^2\right].
\end{align*}
Thus,
\begin{align*}
\frac12\sum_{i=1}^2\tilde R_{ii}=&(n-1)\sum_{\alpha=1}^p\left(T^{\alpha}\right)^2+\frac12\sum_{\alpha=1}^p\left[(n-2)T^{\alpha}\sum_{i=1}^2\mathring{h}^{\alpha}_{ii}-\sum_{j=1}^n\sum_{i=1}^2\left(\mathring{h}^{\alpha}_{ij}\right)^2\right],\\
\frac{1}{n-2}\sum_{i=3}^n\tilde R_{ii}=&(n-1)\sum_{\alpha=1}^p\left(T^{\alpha}\right)^2+\frac{1}{n-2}\sum_{\alpha=1}^p\left[(n-2)T^{\alpha}\sum_{i=3}^n\mathring{h}^{\alpha}_{ii}-\sum_{j=1}^n\sum_{i=3}^n\left(\mathring{h}^{\alpha}_{ij}\right)^2\right].
\end{align*}
By assumption,
\begin{align*}
\tilde R_{ii}+\tilde R_{jj}\geq 2D,\quad\forall 1\leq i<j\leq n,
\end{align*}
then
\begin{align*}
\dfrac12\sum_{i=1}^2\tilde R_{ii}\geq D,\quad\frac{1}{n-2}\sum_{i=3}^n\tilde R_{ii}\geq D.
\end{align*}
Now for every $\varepsilon\in(0,1]$, we get
\begin{align*}
D\leq&\frac{\varepsilon}{2}\sum_{i=1}^2\tilde R_{ii}+\frac{1-\varepsilon}{n-2}\sum_{i=3}^n\tilde R_{ii}\\
=&(n-1)\sum_{\alpha=1}^p\left(T^{\alpha}\right)^2+\frac{n\varepsilon-2}{2}\sum_{\alpha=1}^pT^{\alpha}\sum_{i=1}^2\mathring{h}^{\alpha}_{ii}-\frac{\varepsilon}{2}\sum_{\alpha=1}^p\sum_{j=1}^n\sum_{i=1}^2\left(\mathring{h}^{\alpha}_{ij}\right)^2\\
 & -\frac{1-\varepsilon}{n-2}\sum_{\alpha=1}^p\sum_{j=1}^n\sum_{i=3}^n\left(\mathring{h}^{\alpha}_{ij}\right)^2\\
=&\varepsilon\left(\tilde R_{1212}+\tilde R_{1234}\right)+(n-1-\varepsilon)\sum_{\alpha=1}^p\left(T^{\alpha}\right)^2+\frac{(n-2)\varepsilon-2}{2}\sum_{\alpha=1}^pT^{\alpha}\sum_{i=1}^2\mathring{h}^{\alpha}_{ii}\\
&-\frac{\varepsilon}{2}\sum_{\alpha=1}^p\left(\sum_{i=1}^2\mathring{h}^{\alpha}_{ii}\right)^2-\frac{\varepsilon}{2}\sum_{\alpha=1}^p\sum_{j=3}^n\sum_{i=1}^2\left(\mathring{h}^{\alpha}_{ij}\right)^2-\frac{1-\varepsilon}{n-2}\sum_{\alpha=1}^p\sum_{j=1}^n\sum_{i=3}^n\left(\mathring{h}^{\alpha}_{ij}\right)^2\\
&-\varepsilon\sum_{\alpha=1}^p\left(\mathring{h}^{\alpha}_{13}\mathring{h}^{\alpha}_{24}-\mathring{h}^{\alpha}_{14}\mathring{h}^{\alpha}_{23}\right)\\
\leq&\varepsilon\left(\tilde R_{1212}+\tilde R_{1234}\right)+(n-1-\varepsilon)\sum_{\alpha=1}^p\left(T^{\alpha}\right)^2+\frac{(n-2)\varepsilon-2}{2}\sum_{\alpha=1}^pT^{\alpha}\sum_{i=1}^2\mathring{h}^{\alpha}_{ii}\\
&-\left[\frac{\varepsilon}{2}+\frac{1-\varepsilon}{(n-2)^2}\right]\sum_{\alpha=1}^p\left(\sum_{i=1}^2\mathring{h}^{\alpha}_{ii}\right)^2\\
\leq&\varepsilon\left(\tilde R_{1212}+\tilde R_{1234}\right)+\frac{\left((n-4)\varepsilon+2\right)^2n^2}{8\left(2+\left(n^2-4n+2\right)\varepsilon\right)}\sum_{\alpha=1}^p\left(T^{\alpha}\right)^2.
\end{align*}
\hfill Q.E.D.

\begin{lemma}\label{lemma3.7}
Let $B$ and $\tilde R$ be as in Lemma \ref{lem:a-2}. Assume that for every orthonormal four-frames $\{e_1,e_2,e_3,e_4\}$,
\begin{align*}
\sum_{i=1}^4\sum_{j=1}^n\tilde R_{ijij}\geq 4D, 
\end{align*}
then for all orthonormal four-frames $\{e_1,e_2,e_3,e_4\}$,
\begin{align*}
\sum_{i=1}^2\sum_{j=3}^4\tilde R_{ijij}-2\tilde R_{1234}\geq 4D-\frac{1}{2}\sum_{\alpha=1}^p\left(H^{\alpha}\right)^2.
\end{align*}
\end{lemma}

\noindent\textbf{Proof:} As notations in the proof of Lemma \ref{lem:a-3}, we get
\begin{align*}
\frac14\sum_{i=1}^4\tilde R_{ii}=&(n-1)\sum_{\alpha=1}^p\left(T^{\alpha}\right)^2+\frac14\sum_{\alpha=1}^p\left[(n-2)T^{\alpha}\sum_{i=1}^4\mathring{h}^{\alpha}_{ii}-\sum_{j=1}^n\sum_{i=1}^4\left(\mathring{h}^{\alpha}_{ij}\right)^2\right]\\
=&(n-1)\sum_{\alpha=1}^p\left(T^{\alpha}\right)^2+\frac{n-2}{4}\sum_{\alpha=1}^pT^{\alpha}\sum_{i=1}^4\mathring{h}^{\alpha}_{ii}\\
&-\frac14\sum_{\alpha=1}^p\sum_{i,j=1}^2\left(\mathring{h}^{\alpha}_{ij}\right)^2-\frac14\sum_{\alpha=1}^p\sum_{i,j=3}^4\left(\mathring{h}^{\alpha}_{ij}\right)^2-\frac12\sum_{\alpha=1}^p\sum_{i=1}^2\sum_{j=3}^4\left(\mathring{h}^{\alpha}_{ij}\right)^2\\
&-\frac14\sum_{\alpha=1}^p\sum_{j=5}^n\sum_{i=1}^4\left(\mathring{h}^{\alpha}_{ij}\right)^2\\
\leq&(n-1)\sum_{\alpha=1}^p\left(T^{\alpha}\right)^2+\frac{n-2}{4}\sum_{\alpha=1}^pT^{\alpha}\sum_{i=1}^4\mathring{h}^{\alpha}_{ii}\\
&-\frac18\sum_{\alpha=1}^p\left(\sum_{i=1}^2\mathring{h}^{\alpha}_{ii}\right)^2-\frac18\sum_{\alpha=1}^p\left(\sum_{i=3}^4\mathring{h}^{\alpha}_{ii}\right)^2\\
&-\frac14\sum_{\alpha=1}^p\left(\mathring{h}^{\alpha}_{13}+\mathring{h}^{\alpha}_{24}\right)^2-\frac14\sum_{\alpha=1}^p\left(\mathring{h}^{\alpha}_{14}-\mathring{h}^{\alpha}_{23}\right)^2.
\end{align*}
Notice that
\begin{align*}
\sum_{i=1}^2\sum_{j=3}^4\tilde R_{ijij}-2\tilde R_{1234}=&4\sum_{\alpha=1}^p\left(T^{\alpha}\right)^2+2\sum_{\alpha=1}^pT^{\alpha}\sum_{i=1}^4\mathring{h}^{\alpha}_{ii}+\sum_{\alpha=1}^p\left(\sum_{i=1}^2\mathring{h}^{\alpha}_{ii}\right)\left(\sum_{j=3}^4\mathring{h}^{\alpha}_{jj}\right)\\
&-\sum_{\alpha=1}^p\left(\mathring{h}^{\alpha}_{13}+\mathring{h}^{\alpha}_{24}\right)^2-\sum_{\alpha=1}^p\left(\mathring{h}^{\alpha}_{14}-\mathring{h}^{\alpha}_{23}\right)^2.
\end{align*}
We obtain
\begin{align*}
D\leq&\frac14\left[\sum_{i=1}^2\sum_{j=3}^4\tilde R_{ijij}-2\tilde R_{1234}\right]+(n-1)\sum_{\alpha=1}^p\left(T^{\alpha}\right)^2+\frac{n-2}{4}\sum_{\alpha=1}^pT^{\alpha}\sum_{i=1}^4\mathring{h}^{\alpha}_{ii}\\
&-\dfrac18\sum_{\alpha=1}^p\left(\sum_{i=1}^4\mathring{h}^{\alpha}_{ii}\right)^2-\sum_{\alpha=1}^p\left(T^{\alpha}\right)^2-\frac12\sum_{\alpha=1}^pT^{\alpha}\sum_{i=1}^4\mathring{h}^{\alpha}_{ii}\\
=&\frac14\left[\sum_{i=1}^2\sum_{j=3}^4\tilde R_{ijij}-2\tilde R_{1234}\right]+(n-2)\sum_{\alpha=1}^p\left(T^{\alpha}\right)^2+\frac{n-4}{4}\sum_{\alpha=1}^pT^{\alpha}\sum_{i=1}^4\mathring{h}^{\alpha}_{ii}-\frac18\sum_{\alpha=1}^p\left(\sum_{i=1}^4\mathring{h}^{\alpha}_{ii}\right)^2\\
\leq&\frac14\left[\sum_{i=1}^2\sum_{j=3}^4\tilde R_{ijij}-2\tilde R_{1234}\right]+\frac{n^2}{8}\sum_{\alpha=1}^p\left(T^{\alpha}\right)^2.
\end{align*}
\hfill Q.E.D.

\vspace{.2in}

\section{Proof of Theorem A}

\vspace{.1in}

In this section, we will prove the differentiable sphere theorems for submanifolds in K\"ahler manifold. 

\vspace{.1in}

\noindent\textbf{Proof of Theorem A:} 
 By (\ref{e1.4}), we have for any vector field $X$ on $N$ that
 \begin{equation}\label{E3.1}
\tilde{K}_{\min}|X|^4\leq K(X)\leq \tilde{K}_{\max}|X|^4.
\end{equation}
By (\ref{E2.2}) and (\ref{E3.1}), we have for orthonormal pair $(X,Y)$ on $N$
\begin{eqnarray*}
32K(X,Y)
&\leq& 3\tilde{K}_{\max}\left(|X+JY|^4+|X-JY|^4\right)\nonumber\\
&  & -\tilde{K}_{\min}\left(|X+Y|^4+|X-Y|^4+4|X|^4+4|Y|^4)\right)\nonumber\\
&=& 24(1+\langle X, JY\rangle^2)\tilde{K}_{\max}-16\tilde{K}_{\min}.
\end{eqnarray*}
Similarly we have
 \begin{equation*}
32K(X,Y)\geq 24(1+\langle X, JY\rangle^2)\tilde{K}_{\min}-16\tilde{K}_{\max}.
\end{equation*}
Therefore, we have
 \begin{equation}\label{E3.2}
 \frac{3}{4}(1+\langle X, JY\rangle^2)\tilde{K}_{\min}-\frac{1}{2}\tilde{K}_{\max}\leq K(X,Y)\leq \frac{3}{4}(1+\langle X, JY\rangle^2)\tilde{K}_{\max}-\frac{1}{2}\tilde{K}_{\min}.
\end{equation}

By (\ref{E2.3}) and (\ref{E3.1}), we have for any orthonormal four-frames $\{X,Y,Z,W\}$ on $N$
\begin{eqnarray}\label{E3.3}
256K(X,Y,Z,W)
&\leq& \tilde{K}_{\max}\left(|X+Z+JY+JW|^4+|X+Z-JY-JW|^4\right.\nonumber\\
&  & \ \ \ \ \ \ \left.  +|X-Z+JY-JW|^4+|X-Z-JY+JW|^4   \right.\nonumber\\
&  & \ \ \ \ \ \ \left.  +|X+W+JY-JZ|^4+|X+W-JY+JZ|^4   \right.\nonumber\\
&  & \ \ \ \ \ \ \left.  +|X-W+JY+JZ|^4+|X-W-JY-JZ|^4   \right)\nonumber\\
&  & -\tilde{K}_{\min}\left(|X+Z+JY-JW|^4+|X+Z-JY+JW|^4\right.\nonumber\\
&  & \ \ \ \ \ \ \ \ \left.  +|X-Z+JY+JW|^4+|X-Z-JY-JW|^4   \right.\nonumber\\
&  & \ \ \ \ \ \ \ \ \left.  +|X+W+JY+JZ|^4+|X+W-JY-JZ|^4   \right.\nonumber\\
&  & \ \ \ \ \ \ \ \ \left.  +|X-W+JY-JZ|^4+|X-W-JY+JZ|^4   \right)\nonumber\\
&=& \tilde{K}_{\max}\left[128+8(\langle X+Z,JY+JW\rangle^2+\langle X-Z,JY-JW\rangle^2  \right.\nonumber\\
&  & \ \ \ \ \ \  \ \ \ \ \ \ \ \ \ \ \ \left.  +\langle X+W,JY-JZ\rangle^2+\langle X-W,JY+JZ\rangle^2)   \right]\nonumber\\
&  & -\tilde{K}_{\min}\left[128+8(\langle X+Z,JY-JW\rangle^2+\langle X-Z,JY+JW\rangle^2  \right.\nonumber\\
&  & \ \ \ \ \ \  \ \ \ \ \ \ \ \ \ \ \ \ \ \left.  +\langle X+W,JY+JZ\rangle^2+\langle X-W,JY-JZ\rangle^2)   \right].
\end{eqnarray}
Similarly, we have
\begin{eqnarray}\label{E3.4}
256K(X,Y,Z,W)
&\geq& \tilde{K}_{\min}\left[128+8(\langle X+Z,JY+JW\rangle^2+\langle X-Z,JY-JW\rangle^2  \right.\nonumber\\
&  & \ \ \ \ \ \  \ \ \ \ \ \ \ \ \ \ \ \left.  +\langle X+W,JY-JZ\rangle^2+\langle X-W,JY+JZ\rangle^2)   \right]\nonumber\\
&  & -\tilde{K}_{\max}\left[128+8(\langle X+Z,JY-JW\rangle^2+\langle X-Z,JY+JW\rangle^2  \right.\nonumber\\
&  & \ \ \ \ \ \  \ \ \ \ \ \ \ \ \ \ \ \ \ \left.  +\langle X+W,JY+JZ\rangle^2+\langle X-W,JY-JZ\rangle^2)   \right].
\end{eqnarray}

\vspace{.1in}

Next we will show that under our assumption, $M\times{\mathbb R}^2$ has nonnegative isotropic curvature, i.e., (\ref{e-BS}) holds for all orthonormal four-frames $\{e_1,e_2,e_3,e_4\}$ and all $\lambda, \mu\in[-1,1]$. For that purpose, we first extend the four-frame $\{e_1,e_2,e_3,e_4\}$ to be an orthonormal frame $\{e_1,\cdots,e_{2m}\}$ of $N$ such that $\{e_1,\cdots,e_{n}\}$ are tangent to $M$ and $\{e_{n+1},\cdots,e_{2m}\}$ are normal to $M$. The Gauss equation (\ref{e-gauss}) implies that
\begin{equation}\label{E3.5}
\tilde R(X,Y,Z,W):=R(X,Y,Z,W)-K(X,Y,Z,W)
\end{equation}
is an algebraic curvature. Lemma \ref{lem:a-2} implies that
\begin{equation*}
\tilde R_{1212}+\tilde R_{1234}\geq \dfrac12\left[\dfrac{\sum_{\alpha=1}^p\left(H^{\alpha}\right)^2}{n-1}-\sum_{i,j=1}^n\sum_{\alpha=1}^p\left(h_{ij}^{\alpha}\right)^2\right]=\frac{1}{2}\left(\frac{|{\bf H}|^2}{n-1}-|{\bf B}|^2\right).
\end{equation*}
Lemma \ref{lem31} implies that for   every orthonormal four-frames $\{e_1,e_2,e_3,e_4\}$ and every $\lambda, \mu\in[-1,1]$
\begin{equation*}
\tilde R_{1313}+\lambda^2\tilde R_{1414}+\mu^2\tilde R_{2323}+\lambda^2\mu^2\tilde R_{2424}-2\lambda\mu \tilde R_{1234}\geq \frac{(1+\lambda^2)(1+\mu^2)}{2}\left(\frac{|{\bf H}|^2}{n-1}-|{\bf B}|^2\right),
\end{equation*}
i.e.,
\begin{eqnarray}\label{E3.7}
& & R_{1313}+\lambda^2R_{1414}+\mu^2R_{2323}+\lambda^2\mu^2R_{2424}-2\lambda\mu R_{1234}\nonumber\\
&\geq& K_{1313}+\lambda^2K_{1414}+\mu^2K_{2323}+\lambda^2\mu^2K_{2424}-2\lambda\mu K_{1234}\nonumber\\
&  &+\frac{1+\lambda^2+\mu^2+\lambda^2\mu^2}{2}\left(\frac{|{\bf H}|^2}{n-1}-|{\bf B}|^2\right).
\end{eqnarray}
By (\ref{e-scalar}), we have
\begin{equation}\label{EEE}
|{\bf B}|^2-\frac{1}{n-1}|{\bf H}|^2=\sum_{i,j=1}^nK_{ijij}+\frac{n-2}{n-1}|{\bf H}|^2-R_M.
\end{equation}
Putting (\ref{EEE}) into (\ref{E3.7}) yields
\begin{eqnarray}\label{EE3.7}
& & R_{1313}+\lambda^2R_{1414}+\mu^2R_{2323}+\lambda^2\mu^2R_{2424}-2\lambda\mu R_{1234}\nonumber\\
&\geq& K_{1313}+\lambda^2K_{1414}+\mu^2K_{2323}+\lambda^2\mu^2K_{2424}-2\lambda\mu K_{1234}\nonumber\\
&  &+\frac{1+\lambda^2+\mu^2+\lambda^2\mu^2}{2}\left(R_M-\frac{n-2}{n-1}|{\bf H}|^2-\sum_{i,j=1}^nK_{ijij}\right).
\end{eqnarray}
Therefore, it suffices to estimate the terms involving the curvature tensor $K$ on $N$. 
By (\ref{E3.2}), for every $i\neq j$, we have
 \begin{equation}\label{E3.8}
 K_{ijij}\geq \frac{3}{4}(1+\langle e_i, Je_j\rangle^2)\tilde{K}_{\min}-\frac{1}{2}\tilde{K}_{\max},
\end{equation}
and
 \begin{equation}\label{EE3.8}
 K_{ijij}\leq \frac{3}{4}(1+\langle e_i, Je_j\rangle^2)\tilde{K}_{\max}-\frac{1}{2}\tilde{K}_{\min}.
\end{equation}
Therefore,
\begin{eqnarray}\label{E3.9}
& & K_{1313}+\lambda^2K_{1414}+\mu^2K_{2323}+\lambda^2\mu^2K_{2424}-\frac{1+\lambda^2+\mu^2+\lambda^2\mu^2}{2}\sum_{i,j=1}^nK_{ijij}\nonumber\\
&\geq& (1+\lambda^2+\mu^2+\lambda^2\mu^2)\left(\frac{3}{4}\tilde{K}_{\min}-\frac{1}{2}\tilde{K}_{\max}\right)\nonumber\\
& & +\frac{3}{4}\left(\langle e_1, Je_3\rangle^2+\lambda^2\langle e_1, Je_4\rangle^2+\mu^2\langle e_2, Je_3\rangle^2+\lambda^2\mu^2\langle e_2, Je_4\rangle^2\right)\tilde{K}_{\min}\nonumber\\
& &-\frac{1+\lambda^2+\mu^2+\lambda^2\mu^2}{2}\left[n(n-1)\left(\frac{3}{4}\tilde{K}_{\max}-\frac{1}{2}\tilde{K}_{\min}\right)+\frac{3}{4}\sum_{i,j=1}^n\langle e_i, Je_j\rangle^2\tilde{K}_{\max} \right]\nonumber\\
&=& (1+\lambda^2+\mu^2+\lambda^2\mu^2)\left(\frac{n^2-n+3}{4}\tilde{K}_{\min}-\frac{3n^2-3n+4}{8}\tilde{K}_{\max} \right)\nonumber\\
& & +\frac{3}{4}\left(\langle e_1, Je_3\rangle^2+\lambda^2\langle e_1, Je_4\rangle^2+\mu^2\langle e_2, Je_3\rangle^2+\lambda^2\mu^2\langle e_2, Je_4\rangle^2\right)\tilde{K}_{\min}\nonumber\\
& & -\frac{3(1+\lambda^2+\mu^2+\lambda^2\mu^2)}{8}\sum_{i,j=1}^n\langle e_i, Je_j\rangle^2\tilde{K}_{\max}.
\end{eqnarray}

We will consider three cases:

\vspace{.1in}

{\bf Case 1:} $\tilde{K}_{\min}\geq0$. In this case, we have from (\ref{E3.9}) that
\begin{eqnarray*}
& & K_{1313}+\lambda^2K_{1414}+\mu^2K_{2323}+\lambda^2\mu^2K_{2424}-\frac{1+\lambda^2+\mu^2+\lambda^2\mu^2}{2}\sum_{i,j=1}^nK_{ijij}\nonumber\\
&\geq& (1+\lambda^2+\mu^2+\lambda^2\mu^2)\left(\frac{n^2-n+3}{4}\tilde{K}_{\min}-\frac{3n^2-3n+4}{8}\tilde{K}_{\max} \right)\nonumber\\
& & -\frac{3n(1+\lambda^2+\mu^2+\lambda^2\mu^2)}{8}\tilde{K}_{\max}\nonumber\\
&=& (1+\lambda^2+\mu^2+\lambda^2\mu^2)\left(\frac{n^2-n+3}{4}\tilde{K}_{\min}-\frac{3n^2+4}{8}\tilde{K}_{\max} \right).
\end{eqnarray*}
By (\ref{E3.3}) and (\ref{E3.4}), we have
 \begin{equation}\label{E-mix}
\frac{1}{2}\tilde{K}_{\min}-\tilde{K}_{\max}\leq K_{1234}\leq \tilde{K}_{\max}-\frac{1}{2}\tilde{K}_{\min}.
\end{equation}
Therefore, we have
\begin{eqnarray}\label{E3.10}
& & K_{1313}+\lambda^2K_{1414}+\mu^2K_{2323}+\lambda^2\mu^2K_{2424}-2\lambda\mu K_{1234}\nonumber\\
& & -\frac{1+\lambda^2+\mu^2+\lambda^2\mu^2}{2}\sum_{i,j=1}^nK_{ijij}\nonumber\\
&\geq& (1+\lambda^2+\mu^2+\lambda^2\mu^2)\left(\frac{n^2-n+3}{4}\tilde{K}_{\min}-\frac{3n^2+4}{8}\tilde{K}_{\max}\right)\nonumber\\
& & -\frac{1+\lambda^2+\mu^2+\lambda^2\mu^2}{2}\left(\tilde{K}_{\max}-\frac{1}{2}\tilde{K}_{\min}\right)\nonumber\\
&=& \frac{1+\lambda^2+\mu^2+\lambda^2\mu^2}{2}\left(\frac{n^2-n+4}{2}\tilde{K}_{\min}-\frac{3n^2+8}{4}\tilde{K}_{\max}\right).
\end{eqnarray}
Putting (\ref{E3.10}) into (\ref{EE3.7}) yields
\begin{eqnarray}\label{E3.11}
& & R_{1313}+\lambda^2R_{1414}+\mu^2R_{2323}+\lambda^2\mu^2R_{2424}-2\lambda\mu R_{1234}\nonumber\\
&\geq& \frac{1+\lambda^2+\mu^2+\lambda^2\mu^2}{2}\left(R_M-\frac{3n^2+8}{4}\tilde{K}_{\max}+\frac{n^2-n+4}{2}\tilde{K}_{\min}-\frac{n-2}{n-1}|{\bf H}|^2\right).
\end{eqnarray}

\vspace{.1in}

{\bf Case 2:} $\tilde{K}_{\min}\leq0\leq \tilde{K}_{\max}$. In this case, we have from (\ref{E3.9}) that
\begin{eqnarray*}
& & K_{1313}+\lambda^2K_{1414}+\mu^2K_{2323}+\lambda^2\mu^2K_{2424}-\frac{1+\lambda^2+\mu^2+\lambda^2\mu^2}{2}\sum_{i,j=1}^nK_{ijij}\nonumber\\
&\geq& (1+\lambda^2+\mu^2+\lambda^2\mu^2)\left(\frac{n^2-n+3}{4}\tilde{K}_{\min}-\frac{3n^2-3n+4}{8}\tilde{K}_{\max} \right)\nonumber\\
& & +\frac{3(1+\lambda^2+\mu^2+\lambda^2\mu^2)}{4}\tilde{K}_{\min}-\frac{3n(1+\lambda^2+\mu^2+\lambda^2\mu^2)}{8}\tilde{K}_{\max}\nonumber\\
&=& (1+\lambda^2+\mu^2+\lambda^2\mu^2)\left(\frac{n^2-n+6}{4}\tilde{K}_{\min}-\frac{3n^2+4}{8}\tilde{K}_{\max} \right).
\end{eqnarray*}
By (\ref{E3.3}) and (\ref{E3.4}), we have
 \begin{equation}\label{E-mix2}
\tilde{K}_{\min}-\tilde{K}_{\max}\leq K_{1234}\leq \tilde{K}_{\max}-\tilde{K}_{\min}.
\end{equation}
Therefore, we have
\begin{eqnarray}\label{E3.12}
& & K_{1313}+\lambda^2K_{1414}+\mu^2K_{2323}+\lambda^2\mu^2K_{2424}-2\lambda\mu K_{1234}\nonumber\\
& & -\frac{1+\lambda^2+\mu^2+\lambda^2\mu^2}{2}\sum_{i,j=1}^nK_{ijij}\nonumber\\
&\geq& (1+\lambda^2+\mu^2+\lambda^2\mu^2)\left(\frac{n^2-n+6}{4}\tilde{K}_{\min}-\frac{3n^2+4}{8}\tilde{K}_{\max}\right)\nonumber\\
& & -\frac{1+\lambda^2+\mu^2+\lambda^2\mu^2}{2}\left(\tilde{K}_{\max}-\tilde{K}_{\min}\right)\nonumber\\
&=& \frac{1+\lambda^2+\mu^2+\lambda^2\mu^2}{2}\left(\frac{n^2-n+8}{2}\tilde{K}_{\min}-\frac{3n^2+8}{4}\tilde{K}_{\max}\right).\end{eqnarray}
Putting (\ref{E3.12}) into (\ref{EE3.7}) yields
\begin{eqnarray}\label{E3.13}
& & R_{1313}+\lambda^2R_{1414}+\mu^2R_{2323}+\lambda^2\mu^2R_{2424}-2\lambda\mu R_{1234}\nonumber\\
&\geq& \frac{1+\lambda^2+\mu^2+\lambda^2\mu^2}{2}\left(R_M-\frac{3n^2+8}{4}\tilde{K}_{\max}+\frac{n^2-n+8}{2}\tilde{K}_{\min}-\frac{n-2}{n-1}|{\bf H}|^2\right).
\end{eqnarray}

\vspace{.1in}

{\bf Case 3:} $\tilde{K}_{\max}\leq0$. In this case, we have from (\ref{E3.9}) that
\begin{eqnarray*}
& & K_{1313}+\lambda^2K_{1414}+\mu^2K_{2323}+\lambda^2\mu^2K_{2424}-\frac{1+\lambda^2+\mu^2+\lambda^2\mu^2}{2}\sum_{i,j=1}^nK_{ijij}\nonumber\\
&\geq& (1+\lambda^2+\mu^2+\lambda^2\mu^2)\left(\frac{n^2-n+3}{4}\tilde{K}_{\min}-\frac{3n^2-3n+4}{8}\tilde{K}_{\max} \right)\nonumber\\
& & +\frac{3(1+\lambda^2+\mu^2+\lambda^2\mu^2)}{4}\tilde{K}_{\min}\nonumber\\
&=& (1+\lambda^2+\mu^2+\lambda^2\mu^2)\left(\frac{n^2-n+6}{4}\tilde{K}_{\min}-\frac{3n^2-3n+4}{8}\tilde{K}_{\max} \right).
\end{eqnarray*}
By (\ref{E3.3}) and (\ref{E3.4}), we have
 \begin{equation}\label{E-mix3}
\tilde{K}_{\min}-\frac{1}{2}\tilde{K}_{\max}\leq K_{1234}\leq \frac{1}{2}\tilde{K}_{\max}-\tilde{K}_{\min}.
\end{equation}
Therefore, we have
\begin{eqnarray}\label{E3.14}
& & K_{1313}+\lambda^2K_{1414}+\mu^2K_{2323}+\lambda^2\mu^2K_{2424}-2\lambda\mu K_{1234}\nonumber\\
&\geq& (1+\lambda^2+\mu^2+\lambda^2\mu^2)\left(\frac{n^2-n+6}{4}\tilde{K}_{\min}-\frac{3n^2-3n+4}{8}\tilde{K}_{\max}\right)\nonumber\\
& & -\frac{1+\lambda^2+\mu^2+\lambda^2\mu^2}{2}\left(\frac{1}{2}\tilde{K}_{\max}-\tilde{K}_{\min}\right)\nonumber\\
&=& \frac{1+\lambda^2+\mu^2+\lambda^2\mu^2}{2}\left(\frac{n^2-n+8}{2}\tilde{K}_{\min}-\frac{3n^2-3n+6}{4}\tilde{K}_{\max}\right).
\end{eqnarray}
Putting (\ref{E3.14}) into (\ref{EE3.7}) yields
\begin{eqnarray}\label{E3.15}
& & R_{1313}+\lambda^2R_{1414}+\mu^2R_{2323}+\lambda^2\mu^2R_{2424}-2\lambda\mu R_{1234}\nonumber\\
&\geq& \frac{1+\lambda^2+\mu^2+\lambda^2\mu^2}{2}\left(R_M-\frac{3(n^2-n+2)}{4}\tilde{K}_{\max}+\frac{n^2-n+8}{2}\tilde{K}_{\min}\right.\nonumber\\
& & \hspace{4.3cm} \left.-\frac{n-2}{n-1}|{\bf H}|^2\right).
\end{eqnarray}

\vspace{.1in}

From (\ref{E3.11}), (\ref{E3.13}) and (\ref{E3.15}), we see that in any case, under our assumption (\ref{e-A}), $M\times{\mathbb R}^2$ always has nonnegative isotropic curvature, i.e.,
\begin{equation*}
R_{1313}+\lambda^2R_{1414}+\mu^2R_{2323}+\lambda^2\mu^2R_{2424}-2\lambda\mu R_{1234}\geq0
\end{equation*} 
for all orthonormal four-frames $\{e_1,e_2,e_3,e_4\}$ and all $\lambda, \mu\in[-1,1]$.

\vspace{.1in}

Next, we will estimate the Ricci curvature on $M$. We will assume that $n\geq3$. By the Gauss equation (\ref{e-gauss}), (\ref{e-diff}), (\ref{e-diff2}), (\ref{EEE}), (\ref{E3.8})  and (\ref{EE3.8}), we have for $i\neq j$
\begin{eqnarray}\label{E3.16}
R_{ijij}
&=& K_{ijij}+\sum_{\alpha=n+1}^{2m}[h^{\alpha}_{ii}h^{\alpha}_{jj}-(h^{\alpha}_{ij})^2]\nonumber\\
&\geq& \frac{1}{2}\left(\frac{3}{2}(1+\langle e_i, Je_j\rangle^2)\tilde{K}_{\min}-\tilde{K}_{\max}+R_M-\sum_{i,j=1}^nK_{ijij}-\frac{n-2}{n-1}|{\bf H}|^2\right)\nonumber\\
&\geq& \frac{1}{2}\left(R_M-\frac{3n^2-3n+4}{4}\tilde{K}_{\max}+\frac{n^2-n+3}{2}\tilde{K}_{\min}-\frac{n-2}{n-1}|{\bf H}|^2\right)\nonumber\\
& & +\frac{1}{2}\left(\frac{3}{2}\langle e_i, Je_j\rangle^2\tilde{K}_{\min}-\frac{3}{4}\sum_{i,j=1}^n\langle e_i, Je_j\rangle^2\tilde{K}_{\max}\right),
\end{eqnarray}
with the first equality holds only if 
\begin{equation}\label{e-equa}
h^{\alpha}_{kl}=0, \ \ for \ all \ k\neq l, \{k,l\}\neq \{i,j\} \ and \ any \ \alpha
\end{equation}
and
\begin{equation}\label{e-equa2}
h^{\alpha}_{kk}=h^{\alpha}_{ii}+h^{\alpha}_{jj}, \ \ for \ all \ k\neq i,j,  \ and \ any \ \alpha.
\end{equation}
We will also consider three cases:

\vspace{.1in}

{\bf Case 1:} $\tilde{K}_{\min}\geq0$. In this case, we have from (\ref{E3.16}) and the assumption (\ref{e-A}) that
\begin{equation*}
R_{ijij}
\geq \frac{1}{2}\left(R_M-\frac{3n^2+4}{4}\tilde{K}_{\max}+\frac{n^2-n+3}{2}\tilde{K}_{\min}-\frac{n-2}{n-1}|{\bf H}|^2\right)
\end{equation*} 
with equality holds only if (\ref{e-equa}) and (\ref{e-equa2}) hold. In particular, we see that for any $1\leq i\leq n$,
\begin{equation*}
Ric_{ii}\geq  \frac{n-1}{2}\left(R_M-\frac{3n^2+4}{4}\tilde{K}_{\max}+\frac{n^2-n+3}{2}\tilde{K}_{\min}-\frac{n-2}{n-1}|{\bf H}|^2\right),
\end{equation*}
with equality holds only if 
\begin{equation*}
h^{\alpha}_{ii}=0, \ h^{\alpha}_{kl}=0, \ \ for \ all \ k\neq l, \  and \ any \ \alpha
\end{equation*}
and
\begin{equation*}
h^{\alpha}_{kk}=h^{\alpha}_{ll}, \ \ for \ all \ k,l\neq i,  \ and \ any \ \alpha,
\end{equation*}
which implies
\begin{equation}\label{e-equa3}
|{\bf B}|^2=\frac{|{\bf H}|^2}{n-1}.
\end{equation}
By assumption (\ref{e-A}), we have
\begin{equation*}
Ric_{ii}\geq  \frac{n-1}{2}\left(\tilde{K}_{\max}-\frac{1}{2}\tilde{K}_{\min}\right).
\end{equation*}

\vspace{.1in}

{\bf Case 2:} $\tilde{K}_{\min}\leq0 \leq\tilde{K}_{\max}$. In this case, similar arguments as above shows that \begin{equation*}
Ric_{ii}\geq  \frac{n-1}{2}\left(\tilde{K}_{\max}-\tilde{K}_{\min}\right),
\end{equation*} 
with equality holds only if (\ref{e-equa3}) holds.

\vspace{.1in}

{\bf Case 3:} $\tilde{K}_{\max}\leq0$. In this case, we have from (\ref{E3.16}) and the assumption (\ref{e-A}) that
\begin{equation*}
Ric_{ii}\geq  \frac{n-1}{2}\left(\frac{1}{2}\tilde{K}_{\max}-\tilde{K}_{\min}\right),
\end{equation*} 
with equality holds only if (\ref{e-equa3}) holds.

\vspace{.1in}

If $\tilde K_{\max}$ and $\tilde K_{\min}$ are not both zero, then we can easily see from above that $Ric_M$ is positive everywhere on $M$. 

If $\tilde K_{\max}=\tilde K_{\min}=0$, then by assumption, $M$ has nonnegative Ricci curvature everywhere and has positive Ricci curvature at least at some point. By Aubin's theorem (Lemma \ref{lemma-Aubin}), $M$ admits a metric with positive Ricci curvature. Now we can finish the proof of the theorem:

\vspace{.1in}

If $n=2$, then by our assumption (\ref{e-A}), we see that $M$ has nonnegative Gauss curvature and has positive Gauss curvature at least at some point. Hence $M$ is diffeomorphic to ${\mathbb S}^2$ or ${\mathbb RP}^2$. In particular, since $M$ is simply connected, $M$ is diffeomorphic to ${\mathbb S}^2$.

\vspace{.1in}

If $n=3$, then from the above argument, $M$ admits a metric with positive Ricci curvature. Therefore, $M$ admits a metric with constant positive sectional curvature by Hamilton's theorem (\cite{Ha}). Hence, $M$ is diffeomorphic to a spherical space form. Since $M$ is simply connected, $M$ is diffeomorphic to ${\mathbb S}^3$.

\vspace{.1in}

If $n\geq4$, then $M\times{\mathbb R}^2$ has nonnegative isotropic curvature. On the other hand, putting $\lambda=\mu=1$ in (\ref{E3.7}) and from the above arguments (by considering three cases), we see that under our assumption (\ref{e-A}), we have
\begin{equation*}
R_{1313}+R_{1414}+R_{2323}+R_{2424}-2R_{1234}\geq0
\end{equation*} 
for all orthonormal four-frames $\{e_1,e_2,e_3,e_4\}$. 

\vspace{.1in}

\noindent \textbf{Claim:} {\em $M$ has nonnegative isotropic curvature and has positive isotropic curvature at some point $x_0$ on $M$.}

\vspace{.1in}

\noindent \textbf{Proof of the claim:}
We will also consider three cases according to the sign of the holomorphic sectional curvature as above.

\vspace{.1in}

If $\tilde{K}_{\min}\geq0$, then we have from (\ref{E3.11}), (\ref{e-diff}) and (\ref{e-diff2}) that
\begin{eqnarray*}
& & R_{1313}+R_{1414}+R_{2323}+R_{2424}-2R_{1234}\nonumber\\
&\geq&2\left(R_M-\frac{3n^2+8}{4}\tilde{K}_{\max}+\frac{n^2-n+4}{2}\tilde{K}_{\min}-\frac{n-2}{n-1}|{\bf H}|^2\right)\nonumber\\
&\geq& 0,
\end{eqnarray*}
with the first equality holds only if $h^{\alpha}_{ij}=0$ for all $1\leq i,j\leq n$. We will show that if $\tilde K_{\max}\neq \tilde K_{\min}$ at some point $p\in M$, then the first equality cannot achieve at $p$. Actually, if the first equality holds at $p$, then we have at $p$ that $R_M=\sum_{i,j=1}^nK_{ijij}$ by (\ref{e-scalar}), since $p$ is a totally geodesic point. Now our assumption reduces to 
\begin{equation*}
\sum_{i,j=1}^nK_{ijij}\geq \frac{3n^2+8}{4}\tilde{K}_{\max}-\frac{n^2-n+4}{2}\tilde{K}_{\min}.
\end{equation*} 
Using (\ref{EE3.8}), we compute
\begin{equation*}
\frac{3n^2+8}{4}\tilde{K}_{\max}-\frac{n^2-n+4}{2}\tilde{K}_{\min}
\leq \sum_{i,j=1}^nK_{ijij} \leq \frac{3n^2}{4}\tilde{K}_{\max}-\frac{n^2-n}{2}\tilde{K}_{\min},
\end{equation*}
which implies that $\tilde K_{\max}=\tilde K_{\min}$, contradicting to our assumption. Therefore, if $\tilde K_{\max}\neq \tilde K_{\min}$ at $p$, them $M$ has positive isotropic curvature at $p$. If $\tilde K_{\max}=\tilde K_{\min}$ at $p$, then $M$ has also positive isotropic curvature at $p$ by assumption. 

The proof of the other two cases are similar and we omit the details here. This completes the proof of the claim.
\hfill Q.E.D.

\vspace{.1in}

\noindent\textbf{Proof of Theorem A (continued):} 
By Lemma \ref{lemma-Se} and the above claim, $M$ admits a metric with positive isotropic curvature, and hence $M$ is homeomorphic to a sphere by Micallef-Moore's theorem (Lemma \ref{lemma-MM}). In particular, $M$ is locally irreducible. Now Brendle-Schoen's theorem (Theorem \ref{thmBS}) applying to $M$ gives us that $M$ is either diffeomorphic to a round sphere ${\mathbb S}^n$, or is a K\"ahler manifold biholomorphic to complex projective space, or is isometric to a compact symmetric space. Since, $M$ admits a metric with positive isotropic curvature, Lemma \ref{lemma-MW} shows that $b_2(M)=0$ if $M$ has even dimension, and hence $M$ cannot be a K\"ahler manifold. Furthermore, Seshadri (\cite{Se}) proved that any locally symmetric metric on $M$ must be of constant sectional curvature. Thus, we have shown that $M$ must be diffeomorphic to a round sphere ${\mathbb S}^n$. This finishes the proof of the theorem.
\hfill Q.E.D.

\vspace{.1in}

From the proof of Theorem A, we can easily see that the assumption of Theorem A can be weaken if the submanifold is totally real, which is given by Corollary \ref{cor-totallyreal}.

\vspace{.1in}

\noindent\textbf{Proof of Corollary \ref{cor-totallyreal}:} 
We choose any orthonormal four-frame $\{e_1,e_2,e_3,e_4\}$ on $M$. Since $M$ is totally real in $N$, we see that $Je_i$ is normal to $TM$ for any $1\leq i\leq 4$. Therefore, we have by (\ref{E3.8}) and (\ref{EE3.8}) that for $1\leq i,j\leq 4$
 \begin{equation}\label{E3.17}
K_{ijij}\geq \frac{3}{4}\tilde{K}_{\min}-\frac{1}{2}\tilde{K}_{\max},
\end{equation}
and
 \begin{equation*}
K_{ijij}\leq \frac{3}{4}\tilde{K}_{\max}-\frac{1}{2}\tilde{K}_{\min},
\end{equation*}
Also by (\ref{E3.3}) and (\ref{E3.4}) we have that
\begin{equation}\label{E3.18}
\frac{1}{2}(\tilde{K}_{\min}-\tilde{K}_{\max})\leq K_{1234} \leq \frac{1}{2}(\tilde{K}_{\max}-\tilde{K}_{\min}).
\end{equation}
From (\ref{EE3.7}), (\ref{E3.17}) and (\ref{E3.18}), we see that
\begin{eqnarray*}
& & R_{1313}+\lambda^2R_{1414}+\mu^2R_{2323}+\lambda^2\mu^2R_{2424}-2\lambda\mu R_{1234}\nonumber\\
& \geq& \frac{1+\lambda^2+\mu^2+\lambda^2\mu^2}{2}\left(R_M-\frac{3(n^2-n+2)}{4}\tilde{K}_{\max}+\frac{n^2-n+4}{2}\tilde{K}_{\min}-\frac{n-2}{n-1}|{\bf H}|^2\right).
\end{eqnarray*}
The remaining part of the proof is similar to that of the proof of Theorem A and we omit the details. We only  need to notice that in order to show that the isotropic curvature is nonnegative everywhere and positive at some point on $M$, we have
\begin{eqnarray*}
& & R_{1313}+R_{1414}+R_{2323}+R_{2424}-2R_{1234}\nonumber\\
& \geq& 2\left(R_M-\frac{3(n^2-n+2)}{4}\tilde{K}_{\max}+\frac{n^2-n+4}{2}\tilde{K}_{\min}-\frac{n-2}{n-1}|{\bf H}|^2\right)\nonumber\\
&\geq&0,
\end{eqnarray*}
with the first equality holds at $p\in M$ only if $p$ is a totally geodesic point. Then at $p$, we have $R_M=\sum_{i,j=1}^nK_{ijij}$, and our assumption reduces to
\begin{equation*}
\frac{3(n^2-n+2)}{4}\tilde{K}_{\max}-\frac{n^2-n+4}{2}\tilde{K}_{\min}\leq \sum_{i,j=1}^nK_{ijij} \leq n(n-1)\left(\frac{3}{4}\tilde{K}_{\max}-\frac{1}{2}\tilde{K}_{\min}\right),
\end{equation*}
which implies that $\tilde{K}_{\max}\leq \frac{4}{3}\tilde{K}_{\min}$. But at $p$ we also have
\begin{align*}
R_{1313}+R_{1414}+R_{2323}+R_{2424}-2R_{1234}=&K_{1313}+K_{1414}+K_{2323}+K_{2424}-2K_{1234}\\
\geq&4\tilde K_{\min}-3\tilde K_{\max},
\end{align*}
which implies that $\frac{4}{3}\tilde K_{\min}\leq\tilde K_{\max}$ if $R_{1313}+R_{1414}+R_{2323}+R_{2424}-2R_{1234}=0$. Therefore, $\tilde K_{\max}=\tilde K_{\min}$ at $p$.  
This finished the proof of the corollary.
\hfill Q.E.D.

\vspace{.1in}

\noindent\textbf{Proof of Corollary \ref{cor1.2}:} As in the proof of Corollary \ref{cor-totallyreal}, we have
\begin{eqnarray*}
& & R_{1313}+\lambda^2R_{1414}+\mu^2R_{2323}+\lambda^2\mu^2R_{2424}-2\lambda\mu R_{1234}\nonumber\\
& \geq& \frac{1+\lambda^2+\mu^2+\lambda^2\mu^2}{2}\left(R_M-\frac{n^2-n-2}{4}c-\frac{n-2}{n-1}|{\bf H}|^2\right).
\end{eqnarray*}
It suffices to estimate the isotropic curvature of $M$. By taking $\lambda=\mu=1$, we obtain
\begin{equation*}
R_{1313}+R_{1414}+R_{2323}+R_{2424}-2R_{1234}\geq 2\left(R_M-\frac{n^2-n-2}{4}c-\frac{n-2}{n-1}|{\bf H}|^2\right)\geq0,
\end{equation*}
with the first equality holds at $p\in M$ only if $p$ is a totally geodesic point. Then at $p$, we have $R_M=\sum_{i,j=1}^nK_{ijij}=\frac{n(n-1)c}{4}$. We conclude that $c\geq0$. However,
at $p$, 
\begin{equation*}
R_{1313}+R_{1414}+R_{2323}+R_{2424}-2R_{1234}=K_{1313}+K_{1414}+K_{2323}+K_{2424}-2K_{1234}\geq c.
\end{equation*}
Therefore, if $c\neq 0$, then the isotropic curvature of $M$ is positive everywhere. If $c=0$, then by assumption $M$ has nonnegative isotropic curvature and has positive isotropic curvature at some point $x_0$ on $M$. The remaining part of the proof is similar to that of Theorem A.
\hfill Q.E.D.

\vspace{.2in}

\section{Proof of Theorem B}

\vspace{.1in}

In this section, we will consider differentiable sphere theorem for compact submanifolds in K\"ahler manifold under the Ricci curvature pinching condition.

\vspace{.1in}

\noindent \textbf{Proof of Theorem B:} We will show that under our assumption, $M\times{\mathbb R}^2$ has nonnegative isotropic curvature, i.e., (\ref{e-BS}) holds for all orthonormal four-frames $\{e_1,e_2,e_3,e_4\}$ and all $\lambda, \mu\in[-1,1]$. As in the proof of Theorem A, we first extend the four-frame $\{e_1,e_2,e_3,e_4\}$ to be an orthonormal frame $\{e_1,\cdots,e_{2m}\}$ of $N$ such that $\{e_1,\cdots,e_{n}\}$ are tangent to $M$ and $\{e_{n+1},\cdots,e_{2m}\}$ are normal to $M$. Define the operator $\tilde R$ by (\ref{E3.5}), which is an algebraic curvature. Then for any $1\leq i<j\leq n$, we have from (\ref{EE3.8}) that 
\begin{eqnarray}\label{E4.1}
\sum_{k=1}^n\tilde R_{ikik}+\sum_{k=1}^n\tilde R_{jkjk}
&=&Ric_{ii}+Ric_{jj}-\sum_{k=1}^nK_{ikik}-\sum_{k=1}^nK_{jkjk}        \nonumber\\
&\geq& Ric^{[2]}_{\min}-(n-1)\left(\frac{3}{2}\tilde{K}_{\max}-\tilde{K}_{\min}\right) \nonumber\\
&  & -\frac{3}{4}\sum_{k=1}^n\left(\langle e_i, Je_k\rangle^2+\langle e_j, Je_k\rangle^2\right)\tilde{K}_{\max}.
\end{eqnarray}
Now we will consider three cases:

\vspace{.1in}

{\bf Case 1:} $\tilde{K}_{\min}\geq0$. In this case, we have from (\ref{E4.1}) that
\begin{eqnarray*}
\sum_{k=1}^n\tilde R_{ikik}+\sum_{k=1}^n\tilde R_{jkjk}
\geq  Ric^{[2]}_{\min}-\frac{3n}{2}\tilde{K}_{\max}+(n-1)\tilde{K}_{\min}.
\end{eqnarray*}
By taking $2D=Ric^{[2]}_{\min}-\frac{3n}{2}\tilde{K}_{\max}+(n-1)\tilde{K}_{\min}$ in Lemma \ref{lem:a-3}, we obtain
for every $0<\varepsilon\leq 1$ and all orthonormal four-frames $\{e_1,e_2,e_3,e_4\}$
\begin{equation*}
\tilde R_{1212}+\tilde R_{1234}\geq\dfrac{1}{2\varepsilon}\left[Ric^{[2]}_{\min}-\frac{3n}{2}\tilde{K}_{\max}+(n-1)\tilde{K}_{\min}-\delta(\varepsilon,n)|{\bf H}|^2\right],
\end{equation*}
where $\delta(\varepsilon,n)=\frac{\left((n-4)\varepsilon+2\right)^2}{4\left(2+\left(n^2-4n+2\right)\varepsilon\right)}$.
Lemma \ref{lem31} implies that for every $\lambda, \mu\in[-1,1]$ and every orthonormal four-frames $\{e_1,e_2,e_3,e_4\}$
\begin{eqnarray*}
&  & \tilde R_{1313}+\lambda^2\tilde R_{1414}+\mu^2\tilde R_{2323}+\lambda^2\mu^2\tilde R_{2424}-2\lambda\mu \tilde R_{1234}\nonumber\\
&\geq &\frac{(1+\lambda^2)(1+\mu^2)}{2\varepsilon}\left[Ric^{[2]}_{\min}-\frac{3n}{2}\tilde{K}_{\max}+(n-1)\tilde{K}_{\min}-\delta(\varepsilon,n)|{\bf H}|^2\right],
\end{eqnarray*}
i.e.,
\begin{eqnarray}\label{E4.2}
& & 2\varepsilon(R_{1313}+\lambda^2R_{1414}+\mu^2R_{2323}+\lambda^2\mu^2R_{2424}-2\lambda\mu R_{1234})\nonumber\\
&\geq& 2\varepsilon(K_{1313}+\lambda^2K_{1414}+\mu^2K_{2323}+\lambda^2\mu^2K_{2424}-2\lambda\mu K_{1234})\nonumber\\
&  &+(1+\lambda^2)(1+\mu^2)\left[Ric^{[2]}_{\min}-\frac{3n}{2}\tilde{K}_{\max}+(n-1)\tilde{K}_{\min}-\delta(\varepsilon,n)|{\bf H}|^2\right].
\end{eqnarray}
Since $\tilde{K}_{\min}\geq0$, we have from (\ref{E3.8}) and (\ref{E-mix}) that
\begin{eqnarray}\label{E4.3}
& & K_{1313}+\lambda^2K_{1414}+\mu^2K_{2323}+\lambda^2\mu^2K_{2424}-2\lambda\mu K_{1234}\nonumber\\
&\geq&  (1+\lambda^2)(1+\mu^2)\left(\frac{3}{4}\tilde{K}_{\min}-\frac{1}{2}\tilde{K}_{\max}\right)-\frac{(1+\lambda^2)(1+\mu^2)}{2}\left(\tilde{K}_{\max}-\frac{1}{2}\tilde{K}_{\min}\right)\nonumber\\
&=& (1+\lambda^2)(1+\mu^2)\left(\tilde{K}_{\min}-\tilde{K}_{\max}\right).
\end{eqnarray}
Inserting (\ref{E4.3}) into (\ref{E4.2}), we have
\begin{eqnarray*}
& & 2\varepsilon(R_{1313}+\lambda^2R_{1414}+\mu^2R_{2323}+\lambda^2\mu^2R_{2424}-2\lambda\mu R_{1234})\nonumber\\
&\geq& (1+\lambda^2)(1+\mu^2)\left[Ric^{[2]}_{\min}-\frac{3n+4\varepsilon}{2}\tilde{K}_{\max}+(n-1+2\varepsilon)\tilde{K}_{\min}-\delta(\varepsilon,n)|{\bf H}|^2\right]\nonumber\\
&\geq&0,
\end{eqnarray*}
the strict inequality holds for some point $x_0\in M$, where the last inequality follows from our assumption (\ref{e-B}). The same argument as in the proof of Theorem A implies that $M$ is diffeomorphic to ${\mathbb S}^n$. 

\vspace{.1in}

{\bf Case 2:} $\tilde{K}_{\min}\leq0\leq \tilde {K}_{\max}$. In this case, following the same argument as Case 1, we also have (\ref{E4.2}). By (\ref{E3.8}) and (\ref{E-mix2}), we have
\begin{eqnarray}\label{E4.4}
& & K_{1313}+\lambda^2K_{1414}+\mu^2K_{2323}+\lambda^2\mu^2K_{2424}-2\lambda\mu K_{1234}\nonumber\\
&\geq&  (1+\lambda^2)(1+\mu^2)\left(\frac{3}{2}\tilde{K}_{\min}-\frac{1}{2}\tilde{K}_{\max}\right)-\frac{(1+\lambda^2)(1+\mu^2)}{2}\left(\tilde{K}_{\max}-\tilde{K}_{\min}\right)\nonumber\\
&=& (1+\lambda^2)(1+\mu^2)\left(2\tilde{K}_{\min}-\tilde{K}_{\max}\right).
\end{eqnarray}
Inserting (\ref{E4.4}) into (\ref{E4.2}), we have
\begin{eqnarray*}
& & 2\varepsilon(R_{1313}+\lambda^2R_{1414}+\mu^2R_{2323}+\lambda^2\mu^2R_{2424}-2\lambda\mu R_{1234})\nonumber\\
&\geq& (1+\lambda^2)(1+\mu^2)\left[Ric^{[2]}_{\min}-\frac{3n+4\varepsilon}{2}\tilde{K}_{\max}+(n-1+4\varepsilon)\tilde{K}_{\min}-\delta(\varepsilon,n)|{\bf H}|^2\right]\nonumber\\
&\geq&0,
\end{eqnarray*}
the strict inequality holds for some point $x_0\in M$, where the last inequality follows from our assumption (\ref{e-B}). The same argument as in the proof of Theorem A implies that $M$ is diffeomorphic to ${\mathbb S}^n$. 

\vspace{.1in}

{\bf Case 3:} $\tilde{K}_{\max}\leq0$. In this case, we have from (\ref{E4.1}) that
\begin{eqnarray*}
\sum_{k=1}^n\tilde R_{ikik}+\sum_{k=1}^n\tilde R_{jkjk}
\geq  Ric^{[2]}_{\min}-\frac{3(n-1)}{2}\tilde{K}_{\max}+(n-1)\tilde{K}_{\min}.
\end{eqnarray*}
By taking $2D=Ric^{[2]}_{\min}-\frac{3(n-1)}{2}\tilde{K}_{\max}+(n-1)\tilde{K}_{\min}$ in Lemma \ref{lem:a-3}, we obtain
for every $0<\varepsilon\leq 1$ and all orthonormal four-frames $\{e_1,e_2,e_3,e_4\}$
\begin{equation*}
\tilde R_{1212}+\tilde R_{1234}\geq\dfrac{1}{2\varepsilon}\left[Ric^{[2]}_{\min}-\frac{3(n-1)}{2}\tilde{K}_{\max}+(n-1)\tilde{K}_{\min}-\delta(\varepsilon,n)|{\bf H}|^2\right].
\end{equation*}
Lemma \ref{lem31} implies that for  every orthonormal four-frames $\{e_1,e_2,e_3,e_4\}$ and every $\lambda, \mu\in[-1,1]$
\begin{eqnarray*}
&  & \tilde R_{1313}+\lambda^2\tilde R_{1414}+\mu^2\tilde R_{2323}+\lambda^2\mu^2\tilde R_{2424}-2\lambda\mu \tilde R_{1234}\nonumber\\
&\geq &\frac{(1+\lambda^2)(1+\mu^2)}{2\varepsilon}\left[Ric^{[2]}_{\min}-\frac{3(n-1)}{2}\tilde{K}_{\max}+(n-1)\tilde{K}_{\min}-\delta(\varepsilon,n)|{\bf H}|^2\right],
\end{eqnarray*}
i.e.,
\begin{eqnarray}\label{E4.5}
& & 2\varepsilon(R_{1313}+\lambda^2R_{1414}+\mu^2R_{2323}+\lambda^2\mu^2R_{2424}-2\lambda\mu R_{1234})\nonumber\\
&\geq& 2\varepsilon(K_{1313}+\lambda^2K_{1414}+\mu^2K_{2323}+\lambda^2\mu^2K_{2424}-2\lambda\mu K_{1234})\nonumber\\
&  &+(1+\lambda^2)(1+\mu^2)\left[Ric^{[2]}_{\min}-\frac{3(n-1)}{2}\tilde{K}_{\max}+(n-1)\tilde{K}_{\min}-\delta(\varepsilon,n)|{\bf H}|^2\right].
\end{eqnarray}
Since $\tilde{K}_{\max}\leq0$, we have from (\ref{E3.8}) and (\ref{E-mix3}) that
\begin{eqnarray}\label{E4.6}
& & K_{1313}+\lambda^2K_{1414}+\mu^2K_{2323}+\lambda^2\mu^2K_{2424}-2\lambda\mu K_{1234}\nonumber\\
&\geq&  (1+\lambda^2)(1+\mu^2)\left(\frac{3}{2}\tilde{K}_{\min}-\frac{1}{2}\tilde{K}_{\max}\right)-\frac{(1+\lambda^2)(1+\mu^2)}{2}\left(\frac{1}{2}\tilde{K}_{\max}-\tilde{K}_{\min}\right)\nonumber\\
&=& (1+\lambda^2)(1+\mu^2)\left(2\tilde{K}_{\min}-\frac{3}{4}\tilde{K}_{\max}\right).
\end{eqnarray}
Inserting (\ref{E4.6}) into (\ref{E4.5}), we have
\begin{eqnarray*}
& & 2\varepsilon(R_{1313}+\lambda^2R_{1414}+\mu^2R_{2323}+\lambda^2\mu^2R_{2424}-2\lambda\mu R_{1234})\nonumber\\
&\geq& (1+\lambda^2)(1+\mu^2)\left[Ric^{[2]}_{\min}-\frac{3(n-1+\varepsilon)}{2}\tilde{K}_{\max}+(n-1+4\varepsilon)\tilde{K}_{\min}-\delta(\varepsilon,n)|{\bf H}|^2\right]\nonumber\\
&\geq&0,
\end{eqnarray*}
the strict inequality holds for some point $x_0\in M$, where the last inequality follows from our assumption (\ref{e-B}). The same argument as in the proof of Theorem A implies that $M$ is diffeomorphic to ${\mathbb S}^n$. 
This finishes the proof of the theorem.
\hfill Q.E.D.

\vspace{.1in}

\noindent \textbf{Proof of Corollary \ref{cor1.7}:} Let $M^n$ be a totally real submanifold of a K\"ahler manifold $N^{2m}$. Using the notations as in the proof of Theorem B, we have from (\ref{E4.1}) that
\begin{eqnarray*}
\sum_{k=1}^n\tilde R_{ikik}+\sum_{k=1}^n\tilde R_{jkjk}
&\geq& Ric^{[2]}_{\min}-(n-1)\left(\frac{3}{2}\tilde{K}_{\max}-\tilde{K}_{\min}\right).
\end{eqnarray*}
By taking $2D=Ric^{[2]}_{\min}-\frac{3(n-1)}{2}\tilde{K}_{\max}+(n-1)\tilde{K}_{\min}$ in Lemma \ref{lem:a-3}, we obtain
for every $0<\varepsilon\leq 1$ and all orthonormal four-frames $\{e_1,e_2,e_3,e_4\}$
\begin{equation*}
\tilde R_{1212}+\tilde R_{1234}\geq\dfrac{1}{2\varepsilon}\left[Ric^{[2]}_{\min}-\frac{3(n-1)}{2}\tilde{K}_{\max}+(n-1)\tilde{K}_{\min}-\delta(\varepsilon,n)|{\bf H}|^2\right].
\end{equation*}
Lemma \ref{lem31} implies that for every orthonormal four-frames $\{e_1,e_2,e_3,e_4\}$ and every $\lambda, \mu\in[-1,1]$
\begin{eqnarray*}
&  & \tilde R_{1313}+\lambda^2\tilde R_{1414}+\mu^2\tilde R_{2323}+\lambda^2\mu^2\tilde R_{2424}-2\lambda\mu \tilde R_{1234}\nonumber\\
&\geq &\frac{(1+\lambda^2)(1+\mu^2)}{2\varepsilon}\left[Ric^{[2]}_{\min}-\frac{3(n-1)}{2}\tilde{K}_{\max}+(n-1)\tilde{K}_{\min}-\delta(\varepsilon,n)|{\bf H}|^2\right],
\end{eqnarray*}
i.e.,
\begin{eqnarray}\label{E4.8}
& & 2\varepsilon(R_{1313}+\lambda^2R_{1414}+\mu^2R_{2323}+\lambda^2\mu^2R_{2424}-2\lambda\mu R_{1234})\nonumber\\
&\geq& 2\varepsilon(K_{1313}+\lambda^2K_{1414}+\mu^2K_{2323}+\lambda^2\mu^2K_{2424}-2\lambda\mu K_{1234})\nonumber\\
&  &+(1+\lambda^2)(1+\mu^2)\left[Ric^{[2]}_{\min}-\frac{3(n-1)}{2}\tilde{K}_{\max}+(n-1)\tilde{K}_{\min}-\delta(\varepsilon,n)|{\bf H}|^2\right].
\end{eqnarray}
By (\ref{E3.17}) and (\ref{E3.18}), we have
\begin{eqnarray}\label{E4.9}
& & K_{1313}+\lambda^2K_{1414}+\mu^2K_{2323}+\lambda^2\mu^2K_{2424}-2\lambda\mu K_{1234}\nonumber\\
&\geq&  (1+\lambda^2)(1+\mu^2)\left(\frac{3}{4}\tilde{K}_{\min}-\frac{1}{2}\tilde{K}_{\max}\right)-\frac{(1+\lambda^2)(1+\mu^2)}{2}\left(\frac{1}{2}\tilde{K}_{\max}-\frac{1}{2}\tilde{K}_{\min}\right)\nonumber\\
&=& (1+\lambda^2)(1+\mu^2)\left(\tilde{K}_{\min}-\frac{3}{4}\tilde{K}_{\max}\right).
\end{eqnarray}
Inserting (\ref{E4.9}) into (\ref{E4.8}), we have
\begin{eqnarray*}
& & 2\varepsilon(R_{1313}+\lambda^2R_{1414}+\mu^2R_{2323}+\lambda^2\mu^2R_{2424}-2\lambda\mu R_{1234})\nonumber\\
&\geq& (1+\lambda^2)(1+\mu^2)\left[Ric^{[2]}_{\min}-\frac{3(n-1+\varepsilon)}{2}\tilde{K}_{\max}+(n-1+2\varepsilon)\tilde{K}_{\min}-\delta(\varepsilon,n)|{\bf H}|^2\right]\nonumber\\
&\geq&0,
\end{eqnarray*}
the strict inequality holds for some point $x_0\in M$, where the last inequality follows from our assumption (\ref{e-B}). The same argument as in the proof of Theorem A implies that $M$ is diffeomorphic to ${\mathbb S}^n$. 
This finishes the proof of the corollary.
\hfill Q.E.D.

\vspace{.2in}

\section{Proof of Theorem C and Theorem D}

\vspace{.1in}

In this section, we will prove the topological sphere theorem for submanifolds in K\"ahler manifold.

\vspace{.1in}

\noindent \textbf{Proof of Theorem C:} As before, we will show that under our assumption, $M\times{\mathbb R}^2$ has nonnegative isotropic curvature. For any orthonormal four-frame $\{e_1,e_2,e_3,e_4\}$, we first extend it to be an orthonormal frame $\{e_1,\cdots,e_{2m}\}$ of $N$ such that $\{e_1,\cdots,e_{n}\}$ are tangent to $M$ and $\{e_{n+1},\cdots,e_{2m}\}$ are normal to $M$. The tensor $\tilde R$ defined by (\ref{E3.5}) is an algebraic curvature. Then (\ref{eA2}) and (\ref{e-homm}), (\ref{e-homm2}) implie that
\begin{equation*}
\sum_{i=1}^2\sum_{j=3}^4\tilde R_{ijij}-2\tilde R_{1234}\geq \frac{\sum_{\alpha=1}^p\left(H^{\alpha}\right)^2}{n-2}-\sum_{i,j=1}^n\sum_{\alpha=1}^p\left(h_{ij}^{\alpha}\right)^2=\frac{|{\bf H}|^2}{n-2}-|{\bf B}|^2.
\end{equation*}
i.e.,
\begin{eqnarray}\label{EE6.1}
R_{1313}+R_{1414}+R_{2323}+R_{2424}-2R_{1234}
&\geq& K_{1313}+K_{1414}+K_{2323}+K_{2424}-2K_{1234}\nonumber\\
&  &+\frac{|{\bf H}|^2}{n-2}-|{\bf B}|^2.
\end{eqnarray}
Putting (\ref{EEE}) into (\ref{EE6.1}) yields
\begin{eqnarray}\label{EE6.2}
 R_{1313}+R_{1414}+R_{2323}+R_{2424}-2R_{1234}
&\geq& K_{1313}+K_{1414}+K_{2323}+K_{2424}-2K_{1234}\nonumber\\
&  &+R_M-\frac{n-3}{n-2}|{\bf H}|^2-\sum_{i,j=1}^nK_{ijij}.
\end{eqnarray}
Therefore, it suffices to estimate the terms involving the curvature tensor $K$ on $N$.  
As in the proof of Theorem A, we will consider three cases:

\vspace{.1in}

{\bf Case 1:} $\tilde{K}_{\min}\geq0$. In this case, we have from (\ref{E3.8}), (\ref{EE3.8}) and (\ref{E-mix}) that
\begin{eqnarray}\label{EE6.4}
& & K_{1313}+K_{1414}+K_{2323}+K_{2424}-2K_{1234}-\sum_{i,j=1}^nK_{ijij}\nonumber\\
&\geq& 4\left(\frac{3}{4}\tilde{K}_{\min}-\frac{1}{2}\tilde{K}_{\max}\right)-2\left(\tilde{K}_{\max}-\frac{1}{2}\tilde{K}_{\min}\right)-n(n-1)\left(\frac{3}{4}\tilde{K}_{\max}-\frac{1}{2}\tilde{K}_{\min}\right)\nonumber\\
& & -\frac{3}{4}\sum_{i,j=1}^n\langle e_i,Je_j\rangle^2\tilde{K}_{\max}\nonumber\\
&\geq& \frac{n^2-n+8}{2}\tilde{K}_{\min}-\frac{3n^2+16}{4}\tilde{K}_{\max}.
\end{eqnarray}
Putting (\ref{EE6.4}) into (\ref{EE6.2}) yields
\begin{eqnarray}\label{EE6.5}
 & & R_{1313}+R_{1414}+R_{2323}+R_{2424}-2R_{1234}\nonumber\\
&\geq& R_M-\frac{3n^2+16}{4}\tilde{K}_{\max}+\frac{n^2-n+8}{2}\tilde{K}_{\min}-\frac{n-3}{n-2}|{\bf H}|^2.
\end{eqnarray}

\vspace{.1in}

{\bf Case 2:} $\tilde{K}_{\min}\leq0\leq \tilde{K}_{\max}$. In this case, we have from (\ref{E3.8}), (\ref{EE3.8}) and (\ref{E-mix2}) that
\begin{eqnarray}\label{EE6.6}
& & K_{1313}+K_{1414}+K_{2323}+K_{2424}-2K_{1234}-\sum_{i,j=1}^nK_{ijij}\nonumber\\
&\geq& 4\left(\frac{3}{2}\tilde{K}_{\min}-\frac{1}{2}\tilde{K}_{\max}\right)-2\left(\tilde{K}_{\max}-\tilde{K}_{\min}\right)-n(n-1)\left(\frac{3}{2}\tilde{K}_{\max}-\frac{1}{2}\tilde{K}_{\min}\right)\nonumber\\
& & -\frac{3}{4}\sum_{i,j=1}^n\langle e_i,Je_j\rangle^2\tilde{K}_{\max}\nonumber\\
&=& \frac{n^2-n+16}{2}\tilde{K}_{\min}-\frac{3n^2+16}{4}\tilde{K}_{\max}.
\end{eqnarray}
Putting (\ref{EE6.6}) into (\ref{EE6.2}) yields
\begin{eqnarray}\label{EE6.7}
 & & R_{1313}+R_{1414}+R_{2323}+R_{2424}-2R_{1234}\nonumber\\
&\geq& R_M-\frac{3n^2+16}{4}\tilde{K}_{\max}+\frac{n^2-n+16}{2}\tilde{K}_{\min}-\frac{n-3}{n-2}|{\bf H}|^2.
\end{eqnarray}

\vspace{.1in}

{\bf Case 3:} $\tilde{K}_{\max}\leq0$.  In this case, we have from (\ref{E3.8}), (\ref{EE3.8}) and (\ref{E-mix3}) that
\begin{eqnarray}\label{EE6.8}
& & K_{1313}+K_{1414}+K_{2323}+K_{2424}-2K_{1234}-\sum_{i,j=1}^nK_{ijij}\nonumber\\
&\geq& 4\left(\frac{3}{2}\tilde{K}_{\min}-\frac{1}{2}\tilde{K}_{\max}\right)-2\left(\frac{1}{2}\tilde{K}_{\max}-\tilde{K}_{\min}\right)-n(n-1)\left(\frac{3}{4}\tilde{K}_{\max}-\frac{1}{2}\tilde{K}_{\min}\right)\nonumber\\
&=& \frac{n^2-n+16}{2}\tilde{K}_{\min}-\frac{3(n^2-n+4)}{4}\tilde{K}_{\max}.
\end{eqnarray}
Putting (\ref{EE6.8}) into (\ref{EE6.2}) yields
\begin{eqnarray}\label{EE6.9}
 & & R_{1313}+R_{1414}+R_{2323}+R_{2424}-2R_{1234}\nonumber\\
&\geq& R_M-\frac{3(n^2-n+4)}{4}\tilde{K}_{\max}+\frac{n^2-n+16}{2}\tilde{K}_{\min}-\frac{n-3}{n-2}|{\bf H}|^2.
\end{eqnarray}

\vspace{.1in}

From (\ref{EE6.5}), (\ref{EE6.7}) and (\ref{EE6.9}), we see that in any case, under our assumption (\ref{e-C}), $M$ always has nonnegative isotropic curvature and has positive isotropic curvature at some point. By Lemma \ref{lemma-Se}, $M$ admits a metric with positive isotropic curvature. Since $M$ is simply connected, $M$ is homeomorphic to ${\mathbb S}^n$ by Lemma \ref{lemma-MM}.
\hfill Q.E.D. 

\vspace{.1in}

\noindent \textbf{Proof of Corollary \ref{Cor1.6}:} Let $M^n$ be a totally real submanifold of a K\"ahler manifold $N^{2m}$. In this case, (\ref{EE6.2}) is still true. 
By (\ref{E3.17}) and (\ref{E3.18}), we have
\begin{align}
& K_{1313}+K_{1414}+K_{2323}+K_{2424}-2K_{1234}-\sum_{i,j=1}^nK_{ijij}\notag\\
\geq& 4\left(\frac{3}{4}\tilde{K}_{\min}-\frac{1}{2}\tilde{K}_{\max}\right)-2\left(\frac{1}{2}\tilde{K}_{\max}-\frac{1}{2}\tilde{K}_{\min}\right)-n(n-1)\left(\frac{3}{4}\tilde{K}_{\max}-\frac{1}{2}\tilde{K}_{\min}\right)\notag\\
=& \frac{n^2-n+8}{2}\tilde{K}_{\min}-\frac{3(n^2-n+4)}{4}\tilde{K}_{\max}.\label{EE6.10}
\end{align}
Inserting (\ref{EE6.10}) into (\ref{EE6.2}), we have
\begin{eqnarray*}
 & & R_{1313}+R_{1414}+R_{2323}+R_{2424}-2R_{1234}\nonumber\\
&\geq& R_M-\frac{3(n^2-n+4)}{4}\tilde{K}_{\max}+\frac{n^2-n+8}{2}\tilde{K}_{\min}-\frac{n-3}{n-2}|{\bf H}|^2\nonumber\\
&\geq& 0, 
\end{eqnarray*}
and the strict inequality holds for some point $x_0\in M$, where the last inequality follows from our assumption. Then the corollary follows from Lemma \ref{lemma-Se} and Lemma \ref{lemma-MM}.
\hfill Q.E.D. 

\vspace{.1in}

\noindent \textbf{Proof of Theorem D:} Using the same notations as in the proof of Theorem B, we have from (\ref{EE3.8}) 
\begin{align}
\sum_{i=1}^4\sum_{j=1}^n\tilde R_{ijij}
=& \sum_{i=1}^4Ric_{ii}-\sum_{i=1}^4\sum_{j=1}^nK_{ijij}       \notag\\
\geq& Ric^{[4]}_{\min}-(n-1)\left(3\tilde{K}_{\max}-2\tilde{K}_{\min}\right)-\frac{3}{4}\sum_{i=1}^4\sum_{j=1}^n\langle e_i, Je_j\rangle^2\tilde{K}_{\max}.\label{E6.1}
\end{align}

Now we will consider three cases:

\vspace{.1in}

{\bf Case 1:} $\tilde{K}_{\min}\geq0$. In this case, we have from (\ref{E6.1}) that
\begin{eqnarray*}
\sum_{k=1}^n\tilde R_{ikik}+\sum_{k=1}^n\tilde R_{jkjk}
\geq  Ric^{[4]}_{\min}-3n\tilde{K}_{\max}+2(n-1)\tilde{K}_{\min}.
\end{eqnarray*}
By taking $4D=Ric^{[4]}_{\min}-3n\tilde{K}_{\max}+2(n-1)\tilde{K}_{\min}$ in Lemma \ref{lemma3.7}, we obtain for all orthonormal four-frames $\{e_1,e_2,e_3,e_4\}$,
\begin{align*}
\sum_{i=1}^2\sum_{j=3}^4\tilde R_{ijij}-2\tilde R_{1234}\geq Ric^{[4]}_{\min}-3n\tilde{K}_{\max}+2(n-1)\tilde{K}_{\min}-\frac{1}{2}|{\bf H}|^2.
\end{align*}
In other word,
\begin{eqnarray}\label{E6.2}
R_{1313}+R_{1414}+R_{2323}+R_{2424}-2R_{1234}
&\geq& K_{1313}+K_{1414}+K_{2323}+K_{2424}-2K_{1234}\nonumber\\
&  &+Ric^{[4]}_{\min}-3n\tilde{K}_{\max}+2(n-1)\tilde{K}_{\min}-\frac{1}{2}|{\bf H}|^2.
\end{eqnarray}
Since $\tilde{K}_{\min}\geq0$, we have from (\ref{E3.8}) and (\ref{E-mix}) that
\begin{eqnarray}\label{E6.3}
& & K_{1313}+K_{1414}+K_{2323}+K_{2424}-2K_{1234}\nonumber\\
&\geq& 4\left(\frac{3}{4}\tilde{K}_{\min}-\frac{1}{2}\tilde{K}_{\max}\right)-2\left(\tilde{K}_{\max}-\frac{1}{2}\tilde{K}_{\min}\right)\nonumber\\
&=& 4\left(\tilde{K}_{\min}-\tilde{K}_{\max}\right).
\end{eqnarray}
Inserting (\ref{E6.3}) into (\ref{E6.2}), we have
\begin{eqnarray*}
&   & R_{1313}+R_{1414}+R_{2323}+R_{2424}-2R_{1234}\nonumber\\
& \geq & Ric^{[4]}_{\min}-(3n+4)\tilde{K}_{\max}+2(n+1)\tilde{K}_{\min}-\frac{1}{2}|{\bf H}|^2\nonumber\\
&\geq & 0
\end{eqnarray*}
the strict inequality holds for some point $x_0\in M$, where the last inequality follows from our assumption (\ref{e-D}). By Lemma \ref{lemma-Se}, $M$ admits a metric with positive isotropic curvature. Since $M$ is simply connected, $M$ is homeomorphic to ${\mathbb S}^n$ by Lemma \ref{lemma-MM}.

\vspace{.1in}

{\bf Case 2:} $\tilde{K}_{\min}\leq0\leq \tilde {K}_{\max}$. In this case, following the same argument as Case 1, we also have (\ref{E6.2}). By (\ref{E3.8}) and (\ref{E-mix2}), we have
\begin{eqnarray}\label{E6.4}
& &K_{1313}+K_{1414}+K_{2323}+K_{2424}-2K_{1234}\nonumber\\
&\geq&  4\left(\frac{3}{2}\tilde{K}_{\min}-\frac{1}{2}\tilde{K}_{\max}\right)-2\left(\tilde{K}_{\max}-\tilde{K}_{\min}\right)\nonumber\\
&=& 4\left(2\tilde{K}_{\min}-\tilde{K}_{\max}\right).
\end{eqnarray}
Inserting (\ref{E6.4}) into (\ref{E6.2}), we have
\begin{eqnarray*}
&   & R_{1313}+R_{1414}+R_{2323}+R_{2424}-2R_{1234}\nonumber\\
& \geq & Ric^{[4]}_{\min}-(3n+4)\tilde{K}_{\max}+2(n+3)\tilde{K}_{\min}-\frac{1}{2}|{\bf H}|^2\nonumber\\
&\geq & 0
\end{eqnarray*}
the strict inequality holds for some point $x_0\in M$, where the last inequality follows from our assumption (\ref{e-D}). Then the theorem follows from Lemma \ref{lemma-Se} and Lemma \ref{lemma-MM}.

\vspace{.1in}

{\bf Case 3:} $\tilde{K}_{\max}\leq0$. In this case, we have from (\ref{E6.1}) that
\begin{eqnarray*}
\sum_{k=1}^n\tilde R_{ikik}+\sum_{k=1}^n\tilde R_{jkjk}
\geq  Ric^{[4]}_{\min}-3(n-1)\tilde{K}_{\max}+2(n-1)\tilde{K}_{\min}.
\end{eqnarray*}
By taking $4D=Ric^{[4]}_{\min}-3(n-1)\tilde{K}_{\max}+2(n-1)\tilde{K}_{\min}$ in Lemma \ref{lemma3.7}, we obtain for all orthonormal four-frames $\{e_1,e_2,e_3,e_4\}$,
\begin{align*}
\sum_{i=1}^2\sum_{j=3}^4\tilde R_{ijij}-2\tilde R_{1234}\geq Ric^{[4]}_{\min}-3(n-1)\tilde{K}_{\max}+2(n-1)\tilde{K}_{\min}-\frac{1}{2}|{\bf H}|^2.
\end{align*}
In other word,
\begin{eqnarray}\label{E6.5}
& & R_{1313}+R_{1414}+R_{2323}+R_{2424}-2R_{1234}\nonumber\\
&\geq& K_{1313}+K_{1414}+K_{2323}+K_{2424}-2K_{1234}\nonumber\\
&  &+Ric^{[4]}_{\min}-3(n-1)\tilde{K}_{\max}+2(n-1)\tilde{K}_{\min}-\frac{1}{2}|{\bf H}|^2.
\end{eqnarray}
Since $\tilde{K}_{\max}\leq0$, we have from (\ref{E3.8}) and (\ref{E-mix3}) that
\begin{eqnarray}\label{E6.6}
K_{1313}+K_{1414}+K_{2323}+K_{2424}-2K_{1234}
&\geq&  4\left(\frac{3}{2}\tilde{K}_{\min}-\frac{1}{2}\tilde{K}_{\max}\right)-2\left(\frac{1}{2}\tilde{K}_{\max}-\tilde{K}_{\min}\right)\nonumber\\
&=& 8\tilde{K}_{\min}-3\tilde{K}_{\max}.
\end{eqnarray}
Inserting (\ref{E6.6}) into (\ref{E6.5}), we have
\begin{eqnarray*}
&   & R_{1313}+R_{1414}+R_{2323}+R_{2424}-2R_{1234}\nonumber\\
& \geq & Ric^{[4]}_{\min}-3n\tilde{K}_{\max}+2(n+3)\tilde{K}_{\min}-\frac{1}{2}|{\bf H}|^2\nonumber\\
&\geq & 0
\end{eqnarray*}
the strict inequality holds for some point $x_0\in M$, where the last inequality follows from our assumption (\ref{e-D}). Then the theorem follows from Lemma \ref{lemma-Se} and Lemma \ref{lemma-MM}.
This finishes the proof of the theorem.
\hfill Q.E.D.

\vspace{.1in}

\noindent \textbf{Proof of Corollary \ref{cor1.12}:} Let $M^n$ be a totally real submanifold of a K\"ahler manifold $N^{2m}$. Using the notations as in the proof of Theorem D, we have from (\ref{E6.1}) that
\begin{eqnarray*}
\sum_{i=1}^4\sum_{j=1}^n\tilde R_{ijij}
&\geq& Ric^{[4]}_{\min}-(n-1)\left(3\tilde{K}_{\max}-2\tilde{K}_{\min}\right).
\end{eqnarray*}
By taking $4D=Ric^{[4]}_{\min}-(n-1)\left(3\tilde{K}_{\max}-2\tilde{K}_{\min}\right)$ in Lemma \ref{lemma3.7}, we obtain
for every orthonormal four-frames $\{e_1,e_2,e_3,e_4\}$
\begin{equation*}
\sum_{i=1}^2\sum_{j=3}^4\tilde R_{ijij}-2\tilde R_{1234}\geq Ric^{[4]}_{\min}-3(n-1)\tilde{K}_{\max}+2(n-1)\tilde{K}_{\min}-\frac{1}{2}|{\bf H}|^2.
\end{equation*}
In other word, (\ref{E6.5}) is true.
By (\ref{E3.17}) and (\ref{E3.18}), we have
\begin{eqnarray}\label{EE4.9}
& & K_{1313}+K_{1414}+K_{2323}+K_{2424}-2K_{1234}\nonumber\\
&\geq& 4\left(\frac{3}{4}\tilde{K}_{\min}-\frac{1}{2}\tilde{K}_{\max}\right)-2\left(\frac{1}{2}\tilde{K}_{\max}-\frac{1}{2}\tilde{K}_{\min}\right)=4\tilde{K}_{\min}-3\tilde{K}_{\max}.
\end{eqnarray}
Inserting (\ref{EE4.9}) into (\ref{E6.5}), we have
\begin{eqnarray*}
&   & R_{1313}+R_{1414}+R_{2323}+R_{2424}-2R_{1234}\nonumber\\
& \geq & Ric^{[4]}_{\min}-3n\tilde{K}_{\max}+2(n+1)\tilde{K}_{\min}-\frac{1}{2}|{\bf H}|^2\nonumber\\
&\geq & 0
\end{eqnarray*}
the strict inequality holds for some point $x_0\in M$, where the last inequality follows from our assumption. Then the corollary follows from Lemma \ref{lemma-Se} and Lemma \ref{lemma-MM}.
\hfill Q.E.D. 

\vspace{.2in}


\end{document}